\documentclass{authorsversion}

\usepackage{float}
\usepackage{color}
\usepackage[hyphens]{url}
\usepackage[active]{srcltx}
\usepackage{subfigure}
\usepackage{units}
\usepackage{setspace}
\usepackage{microtype}

\begin{document}
\NME{1}{6}{00}{28}{00}

\runningheads{E.~N.~Lages, E.~S.~S.~Silveira, D.~T.~Cintra and A.~C.~Frery}
{Adaptive time integration based on displacement history curvature}

\title{An adaptive time integration strategy based on displacement history curvature}
\author{E.~N.~Lages\corrauth\footnotemark[2], E.\ S.\ S.\ Silveira\footnotemark[3], D.\ T.\ Cintra\footnotemark[4], A.\ C.\ Frery\footnotemark[5]}

\address{Laborat\'orio de Computa{\c c}{\~a}o Cient\'ifica e Visualiza\c c\~ao, Centro de Tecnologia, Universidade Federal de Alagoas, Campus A. C. Sim\~oes - Av. Lourival de Melo Mota, s/n, Tabuleiro do Martins - Macei\'o/AL, Brazil}

\corraddr{Eduardo Nobre Lages, Laborat\'orio de Computa{\c c}{\~a}o Cient\'ifica e Visualiza\c c\~ao, Centro de Tecnologia, Universidade Federal de Alagoas, Campus A. C. Sim\~oes - Av. Lourival de Melo Mota, s/n, Tabuleiro do Martins - Macei\'o/AL, Brazil}

\footnotetext[2]{E-mail: enl@lccv.ufal.br}
\footnotetext[3]{E-mail: eduardosetton@lccv.ufal.br}
\footnotetext[4]{E-mail: diogotc@lccv.ufal.br}
\footnotetext[5]{E-mail: acfrery@lccv.ufal.br}
%\cgsn{Publishing Arts Research Council}{98--1846389}
\received{12 July 2010}
\revised{}
\noaccepted{}

\begin{abstract}
This work introduces a time-adaptive strategy that uses a refinement estimator based on the first Frenet curvature.
In dynamics, a time-adaptive strategy is a mechanism that interactively proposes changes to the time step used in iterative methods of solution. 
These changes aim to improve the relation between quality of response and computational cost.
The method here proposed is suitable for a variety of numerical time integration problems, e.g., in the study of bodies subjected to dynamical loads.
The motion equation in its space-discrete form is used as reference to derive the formulation presented in this paper.
Our method is contrasted with other ones based on local error estimator and apparent frequencies.
We check the performance of our proposal when employed with the central difference, the explicit generalized-$\alpha$ and the Chung-Lee integration methods.
The proposed refinement estimator demands low computational resources, being easily applied to several direct integration methods.
\end{abstract}


\section{INTRODUCTION}

Computational mechanics is a widely used tool to study many problems in engineering and science.
Many numerical methods have been, and are still being, developed to allow computational simulation ranging from simple mechanical mass-spring models to large scale computational intensive problems, as in computational fluid dynamics.
Finite element, discrete element, and boundary element methods are examples of these approaches.
These techniques are usually computer demanding, mostly when used in sophisticated or realistic models.

The computational effort, even on modern computers, is still an issue that can put limits to the kind of simulations that can be performed in acceptable time windows.
As such, numerical integration methods play a central role when searching for faster and more realistic simulations of all kinds of mechanical models.

Direct integration methods allow the study of time dynamics by a discretization of the solution at time intervals, or increments.
On the one hand, the magnitude of these intervals is related to the convergence of the method being, in general, the smaller the better.
On the other hand, the magnitude has a strong impact on both processing time and precision of the obtained solution.
Regular time intervals, i.e., with constant values along the analysis, are frequently used.
Choosing the proper interval has often been done by previous experience~\cite{BERGAN1985}.
If there is no automatic time stepping facility, one has to rely on intuition to suspend the computer run at a certain time, assess the error and then change the step size before the solution process is resumed~\cite{ZIENKIEWICZ1991}.

Adaptivity is a strategy which allows the control of these time intervals during the simulation to improve the relation between the quality of the solution and the required computational time.
To this end, one can use refinement estimators, measures of the distance between numerical solution obtained by the integration process and the exact solution to the problem.

Time-adaptive strategies have been used as tools for optimizing both performance and quality in numerical analysis of dynamical problems.
In~\cite{PARK1980}, the classic Central Difference Method~\cite{Newmark1959} was used in a strategy of adaptation based on the calculation of the largest apparent frequency of the system.
This frequency was measured by the displacement vectors and incremental accelerations to all degrees of freedom, and from its value the time increment to be adopted was defined.
A similar strategy is found in~\cite{Zhang2001}: a refinement estimator was established from incremental information of displacements and velocities.
With this estimator, obtained at each integration step, the size of the time increment, $\Delta t$, is controlled. 
Also using the largest apparent frequencies, a scheme for the automatic calculation on time increment was proposed in~\cite{BERGAN1985}.
A control function was used to avoid constant changing of $\Delta t$.

Signal processing mechanisms and digital filters have been applied to the problem of numerical integration in dynamical systems.
The technique introduced in~\cite{Soderlind2003} was conceived to produce regularized sequences on time increment.
The use of such sequences had a positive impact on the numerical stability, without further increasing the execution time.

Even if the concept of apparent frequency was effective in the analysis of vibration problems, it can not be applied in problems where this value is close to zero.
Aiming to analyze this kind of problems, an adaptive scheme of implicit integration was developed in~\cite{ZIENKIEWICZ1991}.
A simple and low computational cost estimator based in approximate local solutions was presented.
The proposed local error estimator considers a linearized approximation for the reference local solution, and uses the Newmark method~\cite{Newmark1959}.
In~\cite{Choi1996} an error estimator was developed from a quadratic function approximation to the local solution. 
The parameters of the function are defined by values of acceleration at three time steps, and the procedure is not affected by the direct integration used. 

Also using the concept of approximate local solution, Hulbert and Jang~\cite{HULBERT1995} propose an \emph{a posteriori} error estimator for generalized-$\alpha$ method~\cite{Chung1993}, while Chung \emph{et al.}~\cite{Chung2003} develop an equivalent \emph{a priori} estimator for the same method.
A methodology for the formulation of \emph{a posteriori} error estimators is also presented by Romero and Lacoma~\cite{Romero2006}. 
In that work, the concept of local solution is used in the formulation and illustrated by using Newmark scheme. 

An adaptive technique applied to the explicit generalized-$\alpha$ method (EG-$\alpha$)~\cite{Hulbert1996} was presented in~\cite{Gobat2006}.
In that proposal, the time increment is evaluated at each integration step and, in case it is not precise enough to propagate the solution, its value is reduced by a factor of ten.
The process of reduction can be recursive, with a practical limit of four orders of magnitude to the base value of the increment.

This work presents a new proposal for estimating the integration time step, which is based on the first Frenet curvature.
Such proposal can be applied to any direct integration method.
We present numerical evidence that our proposal is more accurate than the methods presented in~\cite{ZIENKIEWICZ1991} and~\cite{PARK1980}.
We also provide a complexity analysis of our algorithm, showing that it is competitive with other approaches.

The rest of this paper unfolds as follows.
Section~\ref{sec:timeadapt} presents the technique of adaptation while Section~\ref{sec:curvature} presents the regularization procedure employed.
In Section~\ref{sec:comparing} we compare the adaptive strategy with existing methods based on Newmark's algorithms family (sections~\ref{subsec:fourwheeldolly} and~\ref{subsec:elastic}), and we verify that its performance does not depend on the integration strategy (Section~\ref{subseq:VaryingIntegrationStrategy}).
The comparisons employ two different problems, namely four-wheel-dolly and elastic collision.
The adaptive algorithm based on curvature and EG-$\alpha$ method are employed in an application involving mooring lines in Section~\ref{sec:mooring}; in order to do so, our proposal was ported to a production system currently in use in the petroleum and gas industry.
The conclusions are presented in Section~\ref{sec:conclu}.

\section{TIME-ADAPTIVE PROPOSAL}\label{sec:timeadapt}

Our proposal can be used to solve solid dynamics described by the motion equation in its space-discrete form:
\begin{equation}
\boldsymbol{M} \boldsymbol{\ddot{d}} + \boldsymbol{f}^{\text{int}} = \boldsymbol{f}^{\text{ext}},
\label{eq:motionequation}
\end{equation}
where $\boldsymbol{M}$ is the mass matrix, $\boldsymbol{\ddot{d}}$ is the acceleration vector, $\boldsymbol{f}^{\text{int}}$ denotes the internal forces, and $\boldsymbol{f}^{\text{ext}}$ is the vector of external forces.

In order to solve Equation~\eqref{eq:motionequation}, the time variable can be discretized by direct integration methods, for instance, seeking for solutions that satisfy it at times 
\begin{equation}
t_{n+1}=t_{n}+ \Delta t (t), 
\label{eq:discretetime}
\end{equation}
where $n$ denotes the number of evaluations of the motion equation.
This discretization is often performed in a regular way (constant $\Delta t$), whereas in this work we propose an adaptive procedure, hence the dependence on $t$, explicit in \eqref{eq:discretetime}.

The variable $\Delta t$ controls how many times the motion equation is solved for a given time interval of analysis. 
In many situations this increment determines the quality of the numerical response.
For explicit methods of numerical time integration, some stability aspects are related with the value of $\Delta t$.
In terms of computational demand of dynamical software, the variable has high impact in the time required to solve numerical simulations.
The determination of the function that describes a time increment history $\Delta t (t)$ is the main objective of a time-adaptive mechanism.

A way to do this is observing the system dynamics and how it is related with the numerical error provided by the time discretization.
Refinement estimators can be defined aiming to quantify the distance between the numerical response obtained and the idealized exact response of the model.
We propose in this paper to use a geometric indicator, namely, first Frenet curvature~\cite{Kuhnel2006}.

\newcommand{\mytau}{\boldsymbol{\tau}}
\newcommand{\mytaupoint}{\dot{\boldsymbol{\tau}}}
\newcommand{\mytautwopoints}{\ddot{\boldsymbol{\tau}}}
The first Frenet curvature is the norm of the unit tangent derivative by the arc length in a curve. 
For a given curve $\mytau$, parametrized by an arbitrary variable $t$, the value of the curvature can be expressed as:
\begin{eqnarray}
k(t)= \sqrt{\frac{ (\mytaupoint \cdotp \mytaupoint)(\mytautwopoints \cdotp \mytautwopoints)-{(\mytaupoint \cdotp \mytautwopoints)}^2   }{{(\mytaupoint \cdotp \mytaupoint)}^{3}}},
\label{eq:curvatureorig}
\end{eqnarray}
where $\mytaupoint$ and $\mytautwopoints$ are the first and the second derivatives in variable $t$, respectively, and ``$\cdotp$'' denotes the dot product.

In our proposal, the curve $\boldsymbol{\tau}$ is parametrized representing the displacement history for a given mechanic system:  
\begin{eqnarray}
\boldsymbol{\tau}(t)= \begin{Bmatrix} t \\ \boldsymbol{d} \end{Bmatrix},
\label{eq:rdetd}
\end{eqnarray}
where the independent variable $t$ represents time and $\boldsymbol{d}$ is the discrete form of the displacement field in generalized coordinates.

The curvature $k(t)$ for the parametrization presented in \eqref{eq:rdetd} is given by:
\newcommand{\dpoint}	{\dot{\boldsymbol{d}}}
\newcommand{\dtwopoints}{\ddot{\boldsymbol{d}}}
\begin{eqnarray}
	k(t)= \sqrt{\frac{ (1 + \dpoint \cdotp \dpoint)(\dtwopoints \cdotp \dtwopoints)-{(\dpoint \cdotp \dtwopoints)}^2   }{{(1+ \dpoint \cdotp \dpoint)}^{3}}},
\label{eq:curvature}
\end{eqnarray}
where $\dpoint$ and $\dtwopoints$ are the velocity and acceleration vectors, respectively.
This expression is stable from the numerical viewpoint in real world cases. 

In systems with a single degree of freedom, expression~\eqref{eq:curvature} boils down to the following relation:
\begin{equation}
k(t)= \frac{
\bigl|\ddot{d}\bigr|
}{
\bigl(1+{\dot{d}}^2\bigr)^{3/2}
},
\label{eq:onedimensionalcase}
\end{equation}
where $\dot{d}$ and $\ddot{d}$ are the velocity and acceleration, respectively.

The information provided by Equation~\eqref{eq:curvature} can be used in order to control the time increment size.
Since the value of curvature is always positive, one can to establish a relation between the two variables using an exponential function:
\begin{eqnarray}
\Delta t=\max\{\Delta t_{\max} \exp\{-b k(t)\}, \Delta t_{\min} \},
\label{eq:deltatxk}
\end{eqnarray}
where $b$ is a positive constant which captures the prior knowledge about the analytic curvature. The bigger the hypothesized curvature, the smaller this constant should be.
Values of $b$ are provided for each of the examples presented in sections~\ref{sec:comparing} and~\ref{sec:mooring}.
Then, for a given value of curvature $k(t)$, the corresponding integration increment varies between a minimum value $\Delta t_{\min}$, in the limit case where the curvature is infinite, and a maximum value $\Delta t_{\max}$, if $k(t)=0$ results from zero acceleration.

There are two extreme situations. 
When the curvature vanishes, the time increment is maximum. 
It indicates that the system dynamic is inertial only, in other words, the studied body has null acceleration.
In an analogy to a curve mapping in a cartesian space, that situation represents the case where the curve $\boldsymbol{\tau}$ is a straight line.
Note that this condition requires lower refinements in the description of the displacement history.
The another extreme situation, where the curvature is infinite, is verified under abrupt kinematics changes of the analyzed system.
The information yields, as the expression suggests, a reduced value of the integration time step.
The values of $\Delta t_{\min}$ and $\Delta t_{\max}$ are left for the user to specify; in our implementation, we used  ${\Delta t}_{\min}={\Delta t}_{crit} / 100$ and ${\Delta t}_{\max} = 0.85 {\Delta t}_{crit}$ with success.

Equation~\eqref{eq:deltatxk} can be coupled to the algorithm responsible for the time integration task.
For any step evaluated by the algorithm, a value of curvature is obtained from the model kinematics and correlated with a time increment.
This procedure automatically adjusts the time integration process to the model dynamics, optimizing both the quality of the numerical response and the processing time.

Note that the curvature is a function that depends only of lower order derivatives of the displacement function.
This fact permits the presented adaptive procedure to be used in several direct integration methods.
Besides, only simple operations are needed to calculate it: three dot product of vectors.
The simplicity in the refinement estimators obtaining is a basic requirement in a time-adaptive process, once it is an additional task to be included in the time integration algorithm~\cite{BERGAN1985}.

\section{CURVATURE REGULARIZATION}\label{sec:curvature}

A noisy curvature function is a frequent issue when multiple degrees of freedom are used to provide a spatial description of a mechanical system.
It occurs because the refinement estimator is scalar and reflects, in a global way, independent movements of several portions of the mechanical system.
A pre-processing of the curvature before its correlation with the integration time increment is advisable in such cases.

The procedure's purpose is the generation of regularized sequences of curvature in order to improve the computational stability.
A computational procedure is stable, or well conditioned, when small changes of parameters have small effects on the computed response~\cite{Soderlind2006}.
Well conditioned algorithms are seldom available in nonlinear dynamics problems that require time-stepping integration; in this scenario, the proposal of stable strategies is desirable. 

A strategy of regularizing the sequence of curvatures is using a filter that captures maximum values in intervals.
For each instant of the time integration, the curvature value is evaluated and then compared to the maximum value obtained so far.
The biggest value defines the time increment $\Delta t$ used in the current integration step, according to Equation~\eqref{eq:deltatxk}.

The time interval of simulation is sub-divided in reference intervals with extension ${\Delta t}_{d\ell}$, where the maximum values are evaluated.
These intervals can be defined taking as reference a maximum time increment $\Delta t_{\max}$ or the increment which ensures convergence $\Delta t_{crit}$.
In our applications, we defined ${\Delta t}_{d\ell}=100\Delta t_{\max}$.

For each sub-interval ${\Delta t}_{d\ell}$, there is a maximum value of curvature.   
Then, for any time $t$, it is possible to estimate the time increment $\Delta t$ considering, not instantaneous values, but representative values of curvature in sub-intervals.  

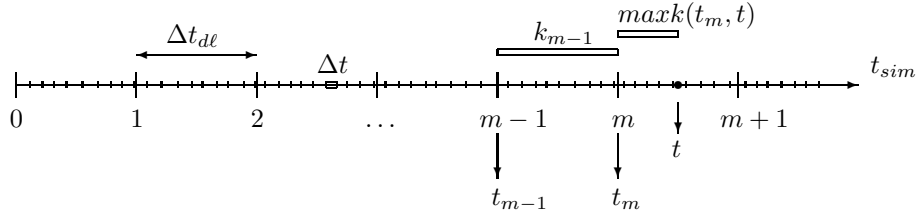
\begin{figure}[hbt]
\centering
\setlength{\unitlength}{0.8cm}
\begin{picture}(16,4)
%timeline
\color[rgb]{0,0,0}
\put(1,2){\vector(1,0){14}}
\put(15.2,2.2){$t_{sim}$}
%{Delta t}_{d\ell} subdivisions
\thinlines
\color[rgb]{0,0,0}
\put(1,2.2){\line(0,-1){0.4}} \put(0.9,1.3){$0$}
\put(3,2.2){\line(0,-1){0.4}} \put(2.9,1.3){$1$}
\put(5,2.2){\line(0,-1){0.4}} \put(4.9,1.3){$2$}
\put(7,2.2){\line(0,-1){0.4}} \put(6.8,1.3){$\dots$}
\put(9,2.2){\line(0,-1){0.4}} \put(8.7,1.3){$m-1$}
\put(11,2.2){\line(0,-1){0.4}} \put(10.9,1.3){$m$}
\put(13,2.2){\line(0,-1){0.4}} \put(12.7,1.3){$m+1$}
%t_m and t_{m-1}
\thinlines
\color[rgb]{0,0,0}
\put(11, 1.2){\vector(0,-1){0.75}}
%\put(10.9,0){$t_m=m {\Delta t}_{d\ell}$}
\put(10.9,0){$t_m$}
\put(9, 1.2){\vector(0,-1){0.75}}
\put(8.9,0){$t_{m-1}$}
%{Delta t}_{d\ell}
\put(4,2.5){\vector(-1,0){1}} \put(4,2.5){\vector(1,0){1}}
\put (3.5,2.7){${\Delta t}_{d\ell}$}
%k_{m-1} e maxk(t_m, t)
\color[rgb]{0,0,0}
\put(9,2.5){\line(0,1){0.1}} \put(9,2.5){\line(1,0){2}}
\put(11,2.6){\line(0,-1){0.1}}\put(11,2.6){\line(-1,0){2}}
\put(9.6,2.75){$k_{m-1}$}
\put(11,2.8){\line(0,1){0.1}} \put(11,2.8){\line(1,0){1}}
\put(12,2.9){\line(0,-1){0.1}}\put(12,2.9){\line(-1,0){1}}
\put(11,3.05){$maxk(t_m,t)$}
%time step intervals (delta_t(t))
\color[rgb]{0,0,0}
\put(1.000000,2.05){\line(0,-1){0.1}} 
\put(1.229167,2.05){\line(0,-1){0.1}} 
\put(1.438228,2.05){\line(0,-1){0.1}} 
\put(1.635105,2.05){\line(0,-1){0.1}} 
\put(1.839860,2.05){\line(0,-1){0.1}} 
\put(2.103529,2.05){\line(0,-1){0.1}} 
\put(2.363434,2.05){\line(0,-1){0.1}} 
\put(2.580029,2.05){\line(0,-1){0.1}} 
\put(2.763672,2.05){\line(0,-1){0.1}} 
\put(2.956250,2.05){\line(0,-1){0.1}} 
\put(3.176331,2.05){\line(0,-1){0.1}} 
\put(3.418100,2.05){\line(0,-1){0.1}} 
\put(3.691843,2.05){\line(0,-1){0.1}} 
\put(3.919886,2.05){\line(0,-1){0.1}} 
\put(4.101098,2.05){\line(0,-1){0.1}} 
\put(4.340837,2.05){\line(0,-1){0.1}} 
\put(4.524061,2.05){\line(0,-1){0.1}} 
\put(4.794194,2.05){\line(0,-1){0.1}} 
\put(4.969252,2.05){\line(0,-1){0.1}} 
\put(5.138656,2.05){\line(0,-1){0.1}} 
\put(5.368610,2.05){\line(0,-1){0.1}} 
\put(5.569021,2.05){\line(0,-1){0.1}} 
\put(5.740588,2.05){\line(0,-1){0.1}} 
\put(5.929999,2.05){\line(0,-1){0.1}} 
\put(6.155469,2.05){\line(0,-1){0.1}} 
\put(6.333889,2.05){\line(0,-1){0.1}} 
\put(6.502868,2.05){\line(0,-1){0.1}} 
\put(6.715140,2.05){\line(0,-1){0.1}} 
\put(6.946380,2.05){\line(0,-1){0.1}} 
\put(7.135204,2.05){\line(0,-1){0.1}} 
\put(7.366680,2.05){\line(0,-1){0.1}} 
\put(7.567790,2.05){\line(0,-1){0.1}} 
\put(7.808382,2.05){\line(0,-1){0.1}} 
\put(7.980656,2.05){\line(0,-1){0.1}} 
\put(8.231766,2.05){\line(0,-1){0.1}} 
\put(8.480437,2.05){\line(0,-1){0.1}} 
\put(8.727048,2.05){\line(0,-1){0.1}} 
\put(9.004378,2.05){\line(0,-1){0.1}} 
\put(9.173838,2.05){\line(0,-1){0.1}} 
\put(9.404960,2.05){\line(0,-1){0.1}} 
\put(9.643073,2.05){\line(0,-1){0.1}} 
\put(9.815146,2.05){\line(0,-1){0.1}} 
\put(10.019520,2.05){\line(0,-1){0.1}} 
\put(10.284773,2.05){\line(0,-1){0.1}} 
\put(10.534320,2.05){\line(0,-1){0.1}} 
\put(10.719817,2.05){\line(0,-1){0.1}} 
\put(10.932252,2.05){\line(0,-1){0.1}} 
\put(11.104902,2.05){\line(0,-1){0.1}} 
\put(11.287444,2.05){\line(0,-1){0.1}} 
\put(11.482120,2.05){\line(0,-1){0.1}} 
\put(11.685131,2.05){\line(0,-1){0.1}} 
\put(11.919574,2.05){\line(0,-1){0.1}} 
\put(12.125360,2.05){\line(0,-1){0.1}} 
\put(12.384298,2.05){\line(0,-1){0.1}} 
\put(12.604080,2.05){\line(0,-1){0.1}} 
\put(12.861618,2.05){\line(0,-1){0.1}} 
\put(13.120751,2.05){\line(0,-1){0.1}} 
\put(13.310843,2.05){\line(0,-1){0.1}} 
\put(13.488288,2.05){\line(0,-1){0.1}} 
\put(13.722604,2.05){\line(0,-1){0.1}} 
\put(13.945783,2.05){\line(0,-1){0.1}} 
\put(14.134771,2.05){\line(0,-1){0.1}} 
\put(14.341550,2.05){\line(0,-1){0.1}} 
%delta t
\put(6.155469,2.05){\line(1,0){0.17842}} 
\put(6.155469,1.95){\line(1,0){0.17842}} 
\put (6.0,2.2){${\Delta t}$}
%t
\thinlines
%\put(12,2.2){\line(0,-1){0.4}}
\put(12,2.0){\circle*{0.15}}
\put(11.9,0.8){$t$}
\put(12,1.7){\vector(0,-1){0.5}}
\end{picture}
\caption{Timeline - Algorithm of regularization by maximum values in intervals.}
\label{fig:timelineMAXKDT}
\end{figure}

Being $m$ the maximum positive integer for which $m {\Delta t}_{d\ell} < t$, and considering the time instants $t_{m-1}=(m-1) {\Delta t}_{d\ell}$ and $t_m=m {\Delta t}_{d\ell}$, it is possible to establish two intervals of reference.
In the first, where $t_{m-1} \leqslant t_{sim}<t_{m}$, we have the maximum value of curvature $k_{m-1}$.
In the second interval, where $t_{m} \leqslant t_{sim} \leqslant t$, the maximum value of curvature is obtained interactively while the time $t$ advances along the timeline $t_{sim}$; the curvature value that represents this last interval is $maxk(t_m,t)$.
Figure~\ref{fig:timelineMAXKDT} illustrates a simulation timeline and the sub-intervals here considered.

The reference to the maximum value of curvature between times $t_{m-1}$ and $t_{m}$ is necessary because the curvature value $maxk(t_m,t)$ has low representativity of the curvature history $k(t)$ when $t$ is close to $t_m$.
A low representativity, or instantaneous values of curvature, may produce frequent changes to the time increment, decreasing, thus, the computational stability.
At the beginning of simulation, the value of $k_{m-1}$ can be adopted null.

Having $k_{m-1}$ and $maxk(t_m,t)$, the value of $k_{m}$ is given by:  
\begin{eqnarray}
k_m=   \Biggl \{ 
	\begin{array}{cl}
	\alpha k_{m-1} + (1- \alpha) maxk(t_m,t)	&\mbox{if $maxk(t_m,t) < k_{m-1}$},\\
	 maxk(t_m,t)					&\mbox{if $maxk(t_m,t) \geqslant k_{m-1}$},\\
	\end{array}
\end{eqnarray}
where $\alpha$ is a weight factor.
A constant $\alpha=1/2$ was employed in all the examples presented in the rest of this article. 
The value of $k_m$ now considers a time regularization procedure and can be used in a correlation with the increment $\Delta t$, according to Equation~\eqref{eq:deltatxk}.  

This regularization procedure requires a negligible computational effort, and produces a stable  and conservative sequence of curvatures.
Additionally, a time step rejection strategy may be used in order to improve the quality of the integration process. 
This can be done by comparing the values of $k_{m-1}$ and $k_{m}$ when $t=t_{m+1}$.
If $k_m > k_{m-1}$, then the time steps between $t_{m}$ and $t_{m+1}$ are discarded and a new integration must be performed for the interval with a time increment $\Delta t (k_m)$.
This procedure skips over the local dynamics of a rapid response change.
Computational codes that implement the strategy must be able to store the entire dynamic system configuration in at least two instants: $t_m$ and $t$.
The change of increment in this case can be done only at instant $t_m$.

The choice of $\zeta$ in the definition of ${\Delta t}_{d\ell}=\zeta\Delta t_{\max}$ has influence on the relative error and on the number of solving steps, and might require a fine tuning; the bigger $\zeta$, the smaller the error and the bigger the number of solving steps.

\begin{figure}[hbt]
\centering
\includegraphics[width=0.7\linewidth]{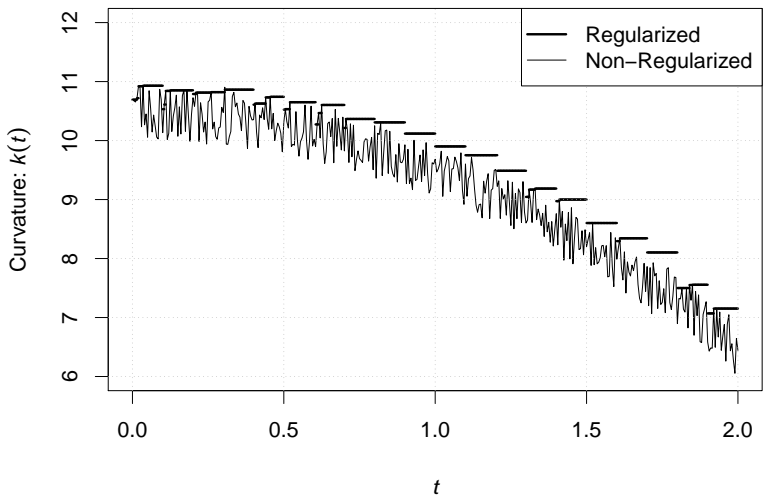}
\caption{Curvature regularization by maximum value in intervals.}
\label{fig:regcurvature}
\end{figure}

Figure~\ref{fig:regcurvature} illustrates the application of the regularization technique to a noisy curvature function.
The hypothetical non-regularized function $k(t)= 10 - t^2 + \beta (t)$ has polynomial basis with a random term $0\leqslant \beta(t) \leqslant 1$.
The regularization interval used is ${\Delta t}_{d\ell} =0.1$, presented for illustration only.

Note that the regularization procedure creates levels within which the time increment does not change.
As commented before, the practice has positive impact in the computational stability.

\section{COMPARISONS}\label{sec:comparing}

This section presents a twofold comparison of the curvature-based strategy with other techniques.
Firstly, sections~\ref{subsec:fourwheeldolly} and~\ref{subsec:elastic} present results fixing the integration strategy (Central Difference Method) and varying the time-adaptive strategy.
Section~\ref{subsec:fourwheeldolly} compares our proposal, the one based on local error and the one based on apparent frequencies for the four-wheel dolly problem.
Section~\ref{subsec:elastic} compares our proposal with local error for the elastic collision of particles; notice that apparent frequencies is inadequate for this problem due to the singularities which appear when the particle is not in contact with the surface.
Secondly, Section~\ref{subseq:VaryingIntegrationStrategy} compares the curvature-based strategy for the same two problems above with three time integration algorithms: the CDM, the EG-$\alpha$, and the Chung-Lee methods.

\subsection{Four-wheel-dolly}\label{subsec:fourwheeldolly}

In this section we apply our proposal to a four-wheel-dolly model excited by an impulsive vertical force.
Underwood and Park~\cite{PARK1980p2} analyzed this problem with an adaptive strategy based on the apparent frequencies of the system~\cite{PARK1980}.
Since their results are not fully reproducible in terms of code and of computational environment~\cite{RRSignalProcessing}, we implemented their strategy in order to compare it with our proposal.

\begin{figure}[hbt]
\centering
\includegraphics[width=0.6\linewidth]{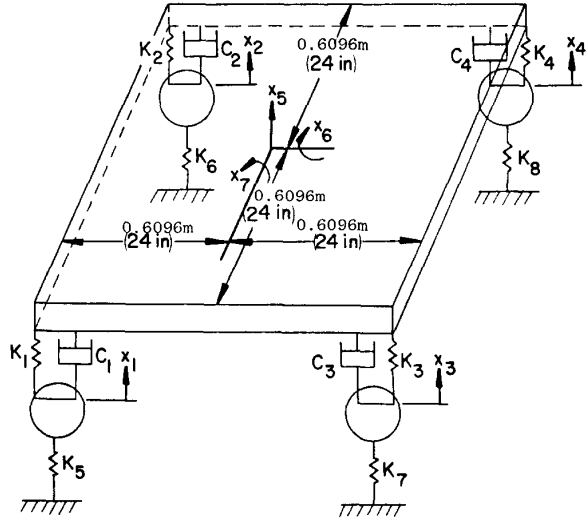}
\caption{Four-wheel-dolly model (As in~\cite{PARK1980p2}).}
\label{fig:modelfwd}
\end{figure}

Figure~\ref{fig:modelfwd} presents the main elements of the problem.
The system has seven degrees of freedom: five translational ($x_1,\dots,x_5$) and two rotational ($x_6,x_7$); each of the four pairs spring-damper ($(k_i,c_i)_{1\leq i\leq 4}$) represents a linear suspension element; $(m_i)_{1\leq i\leq 4}$ represents the four wheels masses; $m_5$ is the translational inertia, while $m_6,m_7$ denote the rotational inertias.
The contact of the dolly with the ground is represented by the four nonlinear springs $(k_i)_{5\leq i\leq 8}$.
The dolly dimensions can be seen in the figure, and the properties are detailed in Table~\ref{tab:propertiesfwd}.

\begin{table}[hbt]
\caption{Dynamic properties of the four-wheel-dolly model (Edited from~\cite{PARK1980p2}).} 
\centering
\begin{tabular}{c r@{.}l rr}
    \toprule
     Number &\multicolumn{2}{c}{Mass} & Dampers & Springs \\
     &\multicolumn{2}{c}{} & (\unit{N s/m}) & (\unit{N/m}) \\
   \midrule                                                                    
    1 & $8$ 	&$7563$~\unit[]{kg} &  700.51	 	 & 87563.43 \\
    2 & $8$ 	&$7563$~\unit[]{kg} &  700.51		 & 87563.43 \\
    3 & $8$ 	&$7563$~\unit[]{kg} &  700.51		 & 87563.43 \\
    4 & $8$ 	&$7563$~\unit[]{kg} &  700.51		 & 87563.43 \\
    5 & $525$ 	&$3804$~\unit[]{kg m${}^2$} & & 175126.85 \\
    6 & $10507$ 	&$6080$~\unit[]{kg m${}^2$} &   & 175126.85 \\
    7 & $10507$ 	&$6080$~\unit[]{kg m${}^2$} &   & 175126.85 \\
    8 &\multicolumn{3}{c}{} 	                   & 175126.85 \\
    \bottomrule
\end{tabular}
\label{tab:propertiesfwd}
\end{table}

The dead load $W$, a constant force of \unit[$-5151.04$]{N} applied to $m_5$, produces the initial displacements $x_1=x_2=x_3=x_4=\unit[-0.007353]{m}$ and $x_5=\unit[-0.22060]{m}$.
The excitation is given by the following expression
\begin{equation}
f_1 = \left\{ 
\begin{array}{ll}
f_{\max}~t/\bar{t} &\text{if~} 0 \leq t \leq \bar{t} \\
f_{\max} \left( 2 - t/\bar{t} \right) & \text{if~} \bar{t} < t \leq 2\bar{t} \\
0 & \text{if~} t > 2\bar{t},
\end{array} 
\right.
\end{equation}
applied to $m_1$, where $f_{\max}=\unit[2224.11]{N}$ and $\bar{t}=\unit[0.025]{s}$.
This impulsive force lifts wheel~1 off the ground and wheel~4 on the rebound.

The general movement description provided by Equation~\eqref{eq:motionequation} is, for the four-wheel-dolly problem under viscous and linear elastic forces, given by:
\begin{equation}
\boldsymbol{M} \boldsymbol{\ddot{d}} + \boldsymbol{C}\boldsymbol{\dot{d}} + \boldsymbol{K}\boldsymbol{d} = \boldsymbol{f}^{\text{ext}},
\end{equation}
where $\boldsymbol{C}$ and $\boldsymbol{K}$ are the damping and stiffness matrices, respectively.

In the following we assume $(k_i)_{1\leq i\leq 4}=k$, $(c_i)_{1\leq i\leq 4}=c$, $(m_i)_{1\leq i\leq 4}=m$, $(x_i=d_i)_{1\leq i\leq 7}$, and for $i=5,6,7,8$:
\begin{equation}
k_i = \left\{ 
\begin{array}{lc}
K &\text{if~} d_{i-4} \leq 0 \\
0 &\text{if~} d_{i-4} > 0.
\end{array} \right.
\end{equation}

Assuming $L=\unit[0.6096]{m}$, one has that
\begin{equation}
\boldsymbol{C}= \left[
\begin{array}{ccccccc}
c & 0 & 0 & 0 & -c & -L c & L c \\
0 & c & 0 & 0 & -c & L c & L c \\
0 & 0 & c & 0 & -c & -L c & -L c \\
0 & 0 & 0 & c & -c & L c & -L c \\
-c & -c & -c & -c & 4c & 0 & 0 \\
-L c & L c & -L c & L c & 0 & 4c L^2 & 0 \\
L c & L c & -L c & -L c & 0 & 0 & 4c L^2
\end{array} 
\right],
\end{equation}
and
\begin{equation}
\boldsymbol{K}= \left[
\begin{array}{ccccccc}
k+k_5 & 0 & 0 & 0 & -k & -L k & L k  \\
0 & k+k_6 & 0 & 0 & -k & L k & L k  \\
0 & 0 & k+k_7 & 0 & -k & -L k & -L k  \\
0 & 0 & 0 & k+k_8 & -k & L k & -L k  \\
-k & -k & -k & -k & 4 k & 0 & 0  \\
-L k & L k & -L k & L k & 0 & 4k L^2 & 0 \\
L k & L k & -L k & -L k & 0 & 0 & 4k L^2
\end{array} 
\right].
\end{equation}
Note that the mass matrix is diagonal, and assuming $m_5=M$ and $m_6=m_7=I$ it reduces to
$\boldsymbol{M}= \lceil 
\begin{array}{ccccccc}
m & m & m & m & M & I & I
\end{array}
\rfloor$, 
and the external force is $\boldsymbol{f}^{\text{ext}} = \left[ 
\begin{array}{ccccccc}
f_1 & 0 & 0 & 0 & -W & 0 & 0
\end{array} 
\right]^{T}$.

Assuming the system is at rest at $t=\unit[0]{s}$, this problem in $\boldsymbol{d}$ is subjected to the initial conditions $\boldsymbol{\dot{d}}=\boldsymbol{0}$ and
\begin{equation}
\boldsymbol{d}= \left[
\begin{array}{ccccccc}
-\frac{W}{4K} & -\frac{W}{4K} & -\frac{W}{4K} & -\frac{W}{4K} & -\frac{W}{4\frac{kK}{k+K}} & 0 & 0
\end{array} \right]^{T}.
\end{equation}

Four integration techniques are compared: the conventional one, with constant time step of comparable length to the other three; an adaptive technique based on natural apparent frequencies~\cite{PARK1980, PARK1980p2}; an adaptive technique based on local errors~\cite{ZIENKIEWICZ1991}; and our proposal which is adaptive and based on the displacement history curvature.
The reference solution was obtained using the conventional solution and a fine time discretization $\Delta t=\unit[10^{-6}]{s}$.

\begin{figure}[hbt]
\centering
\includegraphics[width=0.75\linewidth]{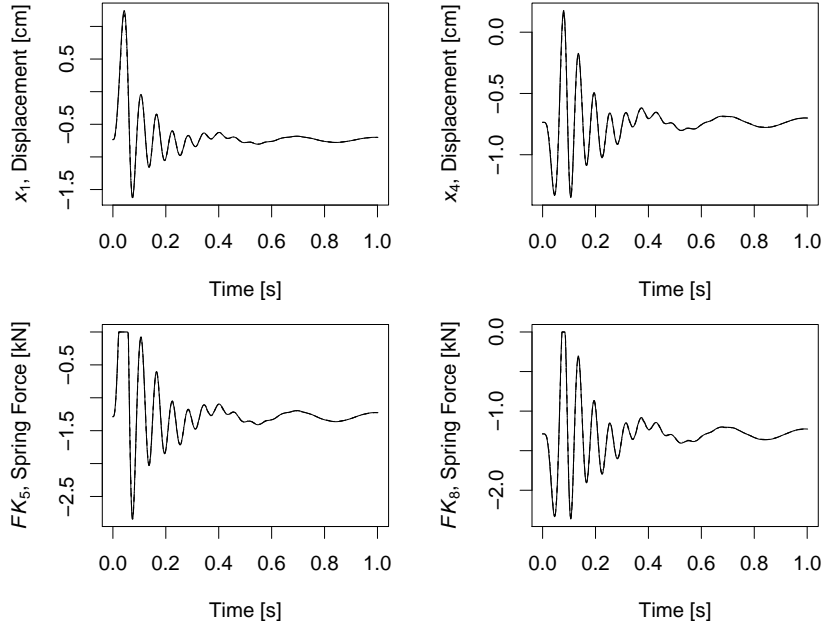}
\caption{Four-wheel-dolly responses.}
\label{fig:fwd_displ_n_forces}
\end{figure}

\begin{figure}[hbt]
\centering
\subfigure[Time step histories.\label{fig:fwd_dthistory}]{\includegraphics[width=0.6\linewidth]{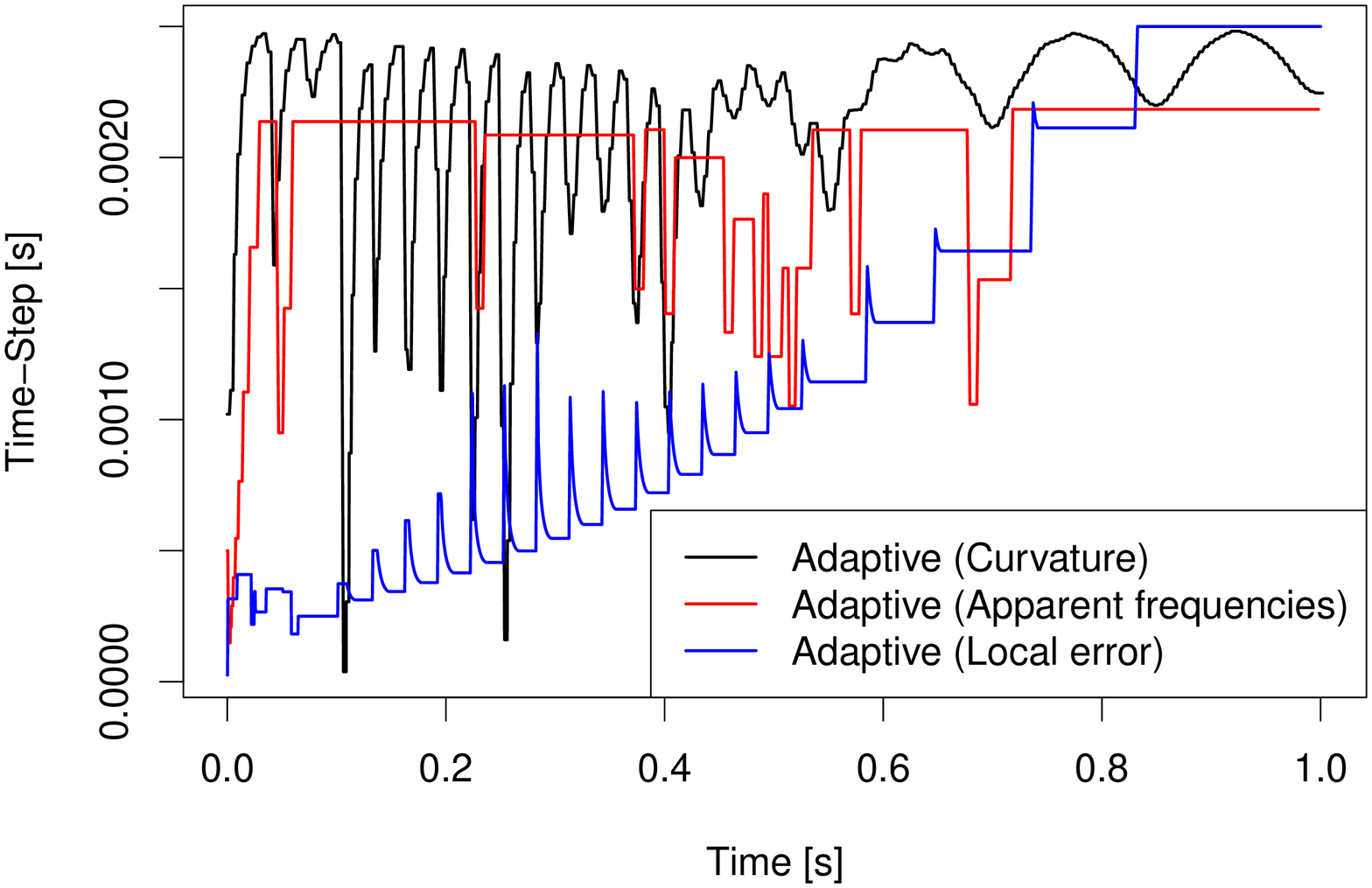}}
\subfigure[Displacement history curvature.\label{fig:fwd_curvature}]{\includegraphics[width=0.6\linewidth]{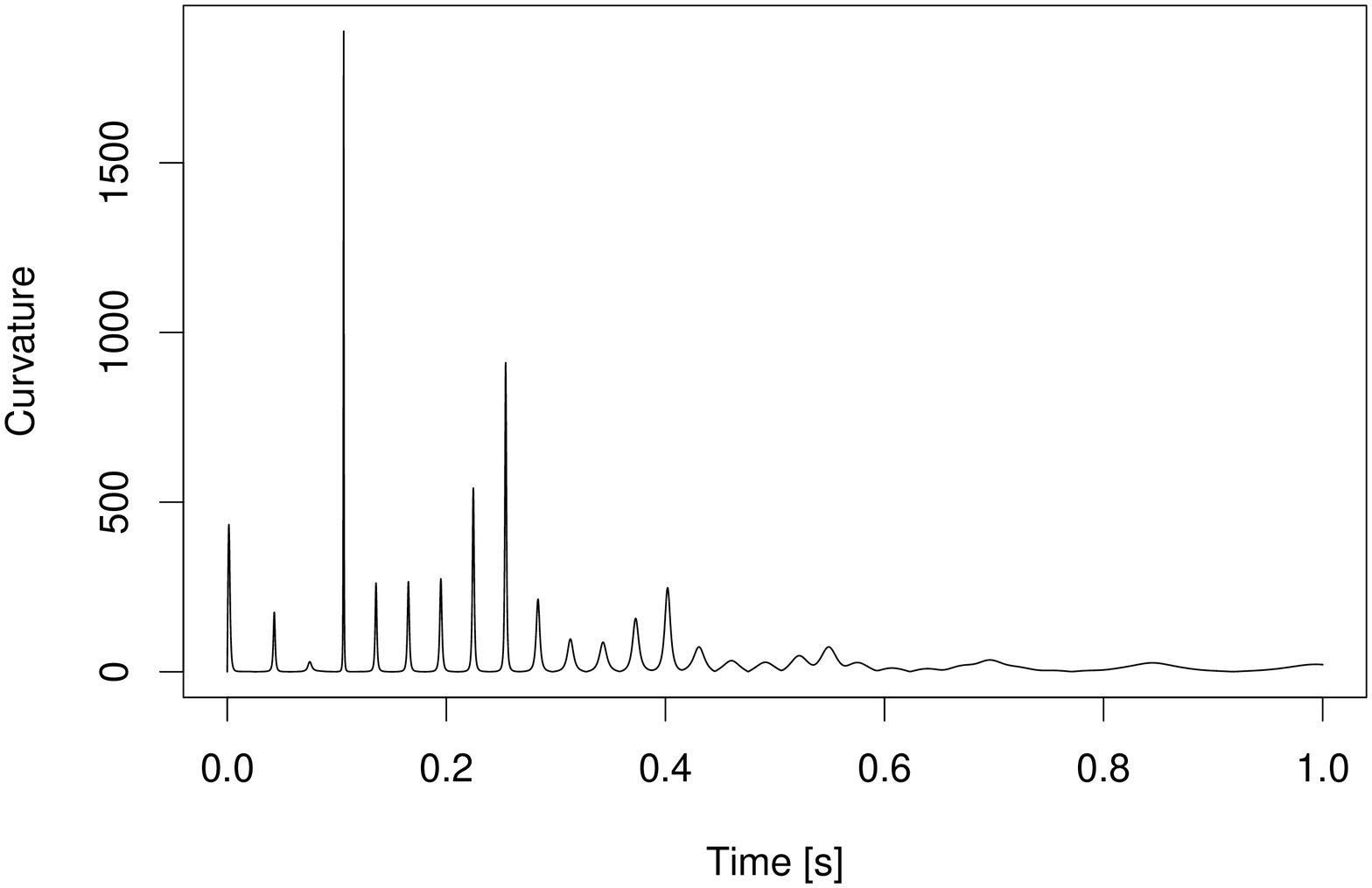}}
\caption{Four-wheel-dolly time step and curvature histories.}
\end{figure}

Figure~\ref{fig:fwd_displ_n_forces} presents the displacements in wheels~1 and~4, and the forces between them and the ground, $FK_5$ and $FK_8$ respectively, as computed by the three adaptive strategies.
Figure~\ref{fig:fwd_dthistory} presents time increment history of both adaptive procedures.
In our proposal, the values $b=0.005$, $\zeta=1$ and $\Delta t_{\max}=\unit[0.0025]{s}$ were used.

The system presents different characteristics along its solution.
In its initial stages, the response is controlled by the highest frequencies that are due to the wheels collisions.
The system tends to obey lower frequencies after the contacts take place and are damped, tending to its original state.
The reference solution of the system shows that wheel~1 (wheel~4 respectively) lifts off the ground at $t=\unit[0.021]{s}$ ($t=\unit[0.073]{s}$, resp.), and rebounds at $t=\unit[0.058]{s}$ ($t=\unit[0.086]{s}$, resp.).
Successful adaptive techniques should be able to detect such changes, and to promote variations in the time step that reduce the number of computations without compromising the quality of the solution.

Figure~\ref{fig:fwd_dthistory} shows the time step histories of the three adaptive algorithms.
It is noticeable that while the apparent frequencies approach exhibits an almost constant behavior with bursts, the algorithm that employs curvature is more closely related to the system nature and to its two main modes: highly nonlinear at the beginning, and progressively more and more linear after the transient.
The algorithm based on local error presents a trend of increasing time-step with bursts which not as closely related to the system behavior as the curvature-based algorithm.
Also, the algorithm based on local error is more time consuming, as will be quantified later, since it takes a relatively long time to increase its time-step.

Figure~\ref{fig:fwd_curvature} presents the displacement history curvature as computed by the reference solution.
Changes in this function, which coincide with those instants when the system is undergoing intense nonlinearities, control the time step of our proposal, as presented in Figure~\ref{fig:fwd_dthistory}.

\begin{figure}[hbt]
\centering
\includegraphics[width=0.75\linewidth]{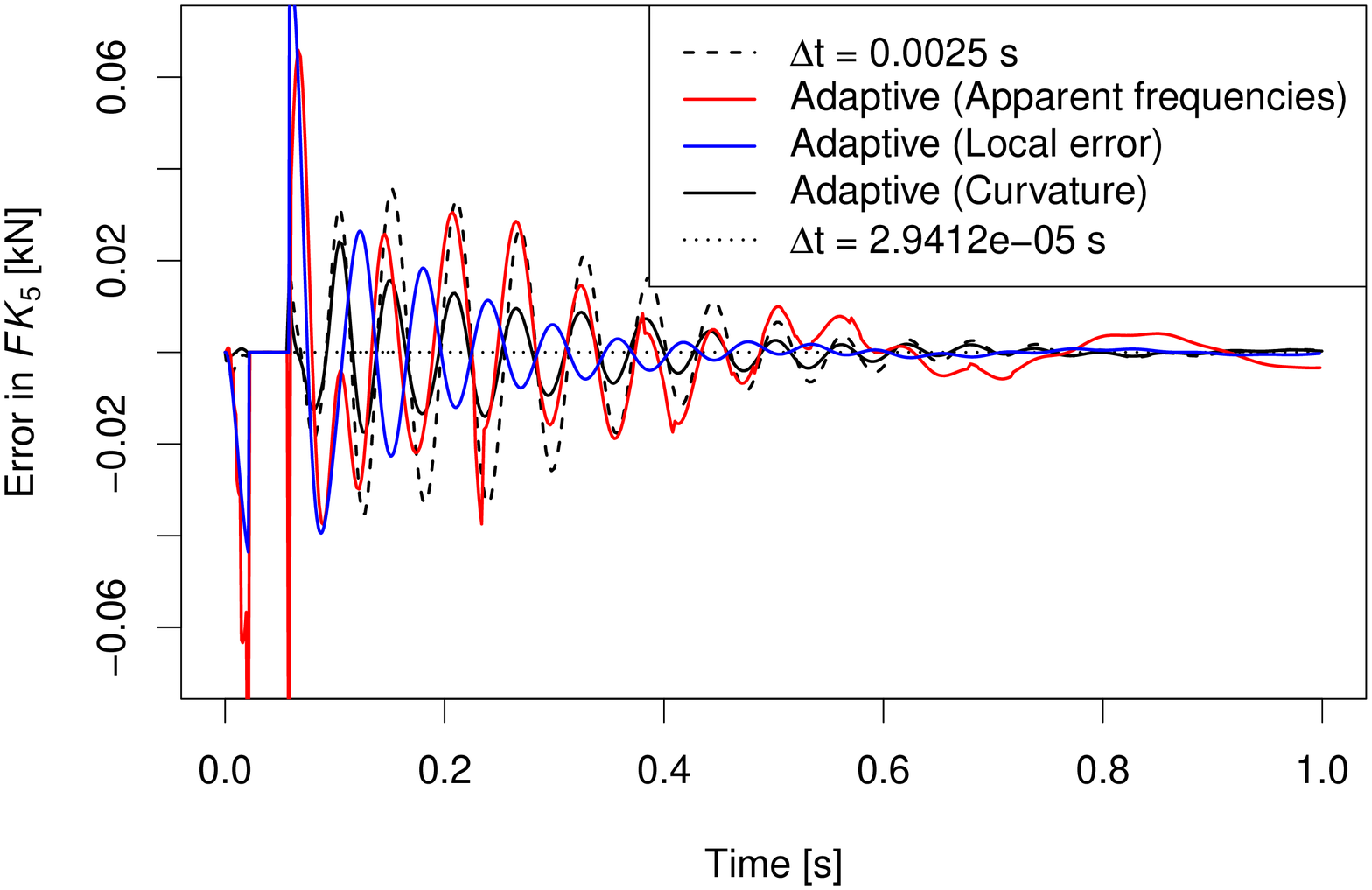}
\caption{Absolute errors computing the force at spring 5.}
\label{fig:fwd_stdvscurv}
\end{figure}

Figure~\ref{fig:fwd_stdvscurv} presents absolute errors of the force in wheel~1.
The errors are measured with respect to the reference solution which was obtained with fixed time step $\Delta t=\unit[10^{-6}]{s}$.
Five errors are shown: three from adaptive strategies and two observed with fixed time steps: $\Delta t_{\min} = \unit[2.9412\cdot10^{-5}]{s}$ and $\Delta t_{\max} = \unit[2.5\cdot10^{-3}]{s}$.
These two last are the smallest and biggest time steps observed in our adaptive solution.

Figure~\ref{fig:fwd_stdvscurv} shows that the smallest error due adaptive algorithms is produced by the proposal based on curvature.
This error is bigger than the one produced by the smallest fixed time step, but smaller than the one yielded by the adaptive algorithms which employ apparent frequencies and local error which, in turn, produce an error smaller than the one associated to the biggest time step.
The relative improvement of the algorithm based on curvature with respect to the one based on apparent frequencies is obtained albeit the former employs time increments which are close to the biggest time steps in a significant portion of the simulation (c.f.\ Figure~\ref{fig:fwd_dthistory}).

\begin{figure}[hbt]
\centering
\includegraphics[width=0.75\linewidth]{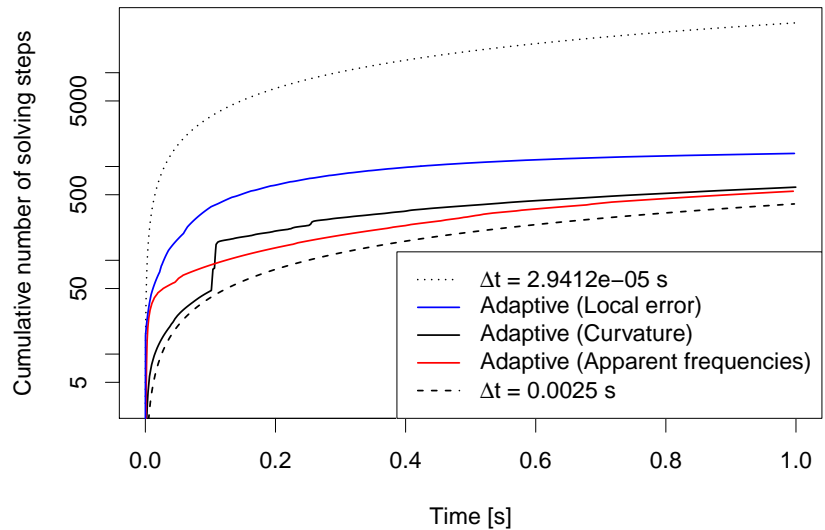}
\caption{Cumulative number of solving steps.}
\label{fig:fwd_solvsteps_hist}
\end{figure}

Figure~\ref{fig:fwd_solvsteps_hist} presents the cumulative number of solving steps of the five strategies.
As expected, the most and least demanding algorithms are, respectively, the ones with fixed time step and smallest and biggest increment.
The adaptive solutions lie between those extremes.
Disregarding those extremes and the initial stages, the cumulative number of solving steps are decreasingly ordered as follows: the one based on local error, the one based on curvature and the one based on apparent frequencies.
The difference between the two last approaches tends to diminish with time, making both approaches comparable in terms of number of solving steps.

\subsection{Elastic Collision of Particles}\label{subsec:elastic}

The following example illustrates the application of the new time-adaptive technique to a problem of elastic collision of particles.
The problem has been inspected by three different methodologies. 
The first uses the analytical solution to investigate the unknown kinematics, and is used as a reference; the other solutions employ the Discrete Element Method~\cite{CUNDALL1979}.

Two different numerical integration procedures are performed by the numerical solutions: the standard, where the increment is constant during entire simulation, and the adaptive, where this increment varies along the simulation based on the displacement curvature history.
In the latter we used the values $b=0.444$, $\zeta=10$, $\Delta t_{crit}=\unit[2.0\cdot10^{-5}]{s}$ and $\Delta t_{\max}=0.85{\Delta t}_{crit}$.
The standard solution employs the Central Difference Method (CDM) to the numerical integration and uses a time step $\Delta t = 0.1{\Delta t}_{crit}$, as defined in Table~\ref{tab:freefallproblem}.

\begin{figure}[hbt]
\centering
\includegraphics[width=0.8\linewidth]{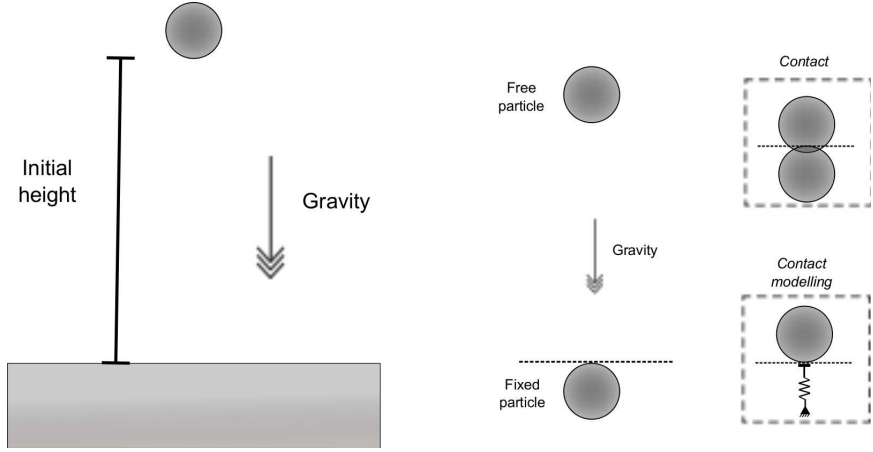}
\caption{Free fall with perfectly elastic collision.}
\label{fig:freefalling}
\end{figure}

The problem is obtaining the motion history of a spherical body subjected to gravitational and contact forces.
The sphere is initially at rest, and suddenly falls from a certain height of a reference surface, as illustrated in Figure~\ref{fig:freefalling}.

The dynamics of this problem is well known.
Due to the gravitational force, the sphere develops a movement until impacting the surface.
After the impact, the body is thrown up and rises up to its initial position.
This movement of rise and fall repeats indefinitely because of the elastic nature of impact.

Being $g$ the gravity acceleration, $h_0$ the initial height and $k$ the elastic stiffness of a spring that represents the contact, the particle vertical position history $h(t)$ can be obtained from the movement differential equation, and is given by
\newcommand{\ufreefall}{h_0- \frac{1}{2}gt^2}
\newcommand{\ucontact}{ \frac{mg}{k} \cos \Bigl( t_c \sqrt{\frac{k}{m}} \Bigr) -\sqrt{ \frac{2mgh}{k}} \sin \Bigl( t_c \sqrt{\frac{k}{m}} \Bigr)  - \frac{mg}{k} }
\newcommand{\urest}{	{t_r}\sqrt{2gh}-\frac{1}{2}g{t_r}^2}
\begin{eqnarray}
h(t) =  \left\{ 
	\begin{array}{cl}
		\ufreefall 			& \mbox{if $t\leq t_q$},\\
		\ucontact 			& \mbox{if $t_q<t \leq t_{ac}$ },\\
		\urest				& \mbox{if $t_{ac} < t \leq t_{f}$},
	\end{array}
	\right.
\label{eq:htparticle}
\end{eqnarray}
where
\begin{align}
t_c &= t - t_q,\\
t_r &= t - t_{ac},\\
t_q = \sqrt{\frac{2h_o}{g}},\quad t_{ac}&= t_q + t_{cont},\quad t_f= 2~t_q + t_{cont},
\end{align}
and $t_{cont}$ is the time interval where there is contact with the surface.
This value is the first positive root of $h(t)=0$ for the second branch of Equation~\eqref{eq:htparticle}.
The repetition period of the movement of rise and fall is given by $t_f$.     

The numerical modeling of the contact with the surface, using discrete elements, is performed assuming the interaction between two distinct particles (Figure~\ref{fig:freefalling}).
Aiming to represent the surface where the contact takes place, and simplifying the contact detection process, one of the particles has fixed position.
The second particle is free to move vertically and its displacement history is investigated.
The physical and geometric parameters used in the numerical analyses are presented in Table~\ref{tab:freefallproblem}.

\begin{table}
\caption{Parameters of the contact numerical model.} 
\begin{center}
\begin{small}
\begin{tabular}{cccc}
    \toprule
    Parameter              & Symbol   & Unit   & Value \\
    \midrule
    Gravity                 & $g$        & $\unitfrac{m}{s^2}$         & $10$ \\
    Launching height         & $h_0$      & $\unit{m}$             &$1.25$\\
    Contact stiffness       & $k_c$      & $\unitfrac{kN}{m}$          & $10^{10}$ \\
    Particle mass	    & $m$        &$\unit{kg}$             &$1$\\
    Critical time step	    & ${\Delta t}_{crit}$        &$\unit{s}$             &$2 \cdot 10^{-5}$\\
    \bottomrule
\end{tabular}
\end{small}
\end{center}
\label{tab:freefallproblem}
\end{table}

\begin{figure}[hbt]
\centering
\subfigure[\label{fig:initial}]{\includegraphics[width=0.45\linewidth]{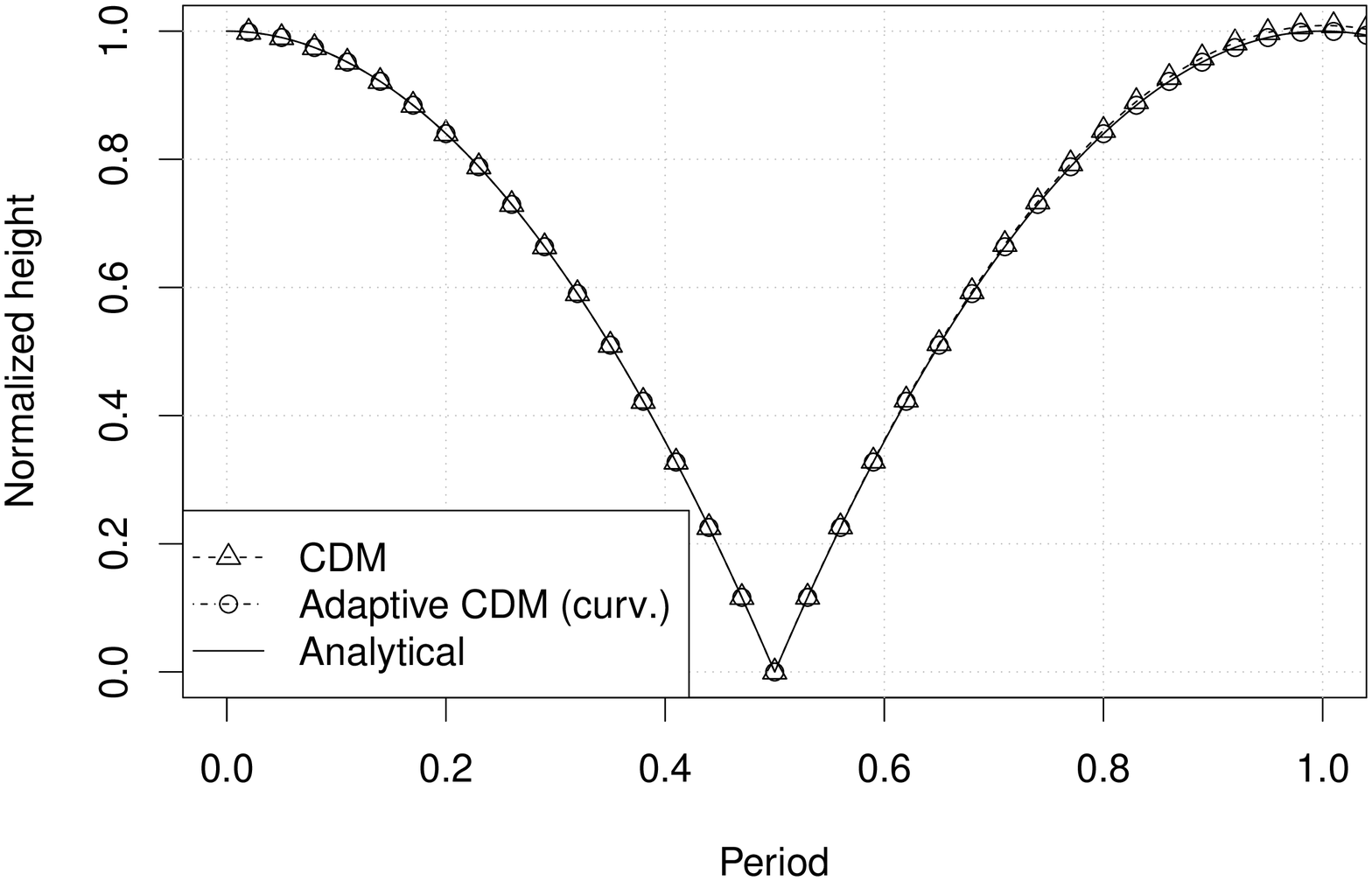}}
\subfigure[\label{fig:interm}]{\includegraphics[width=0.45\linewidth]{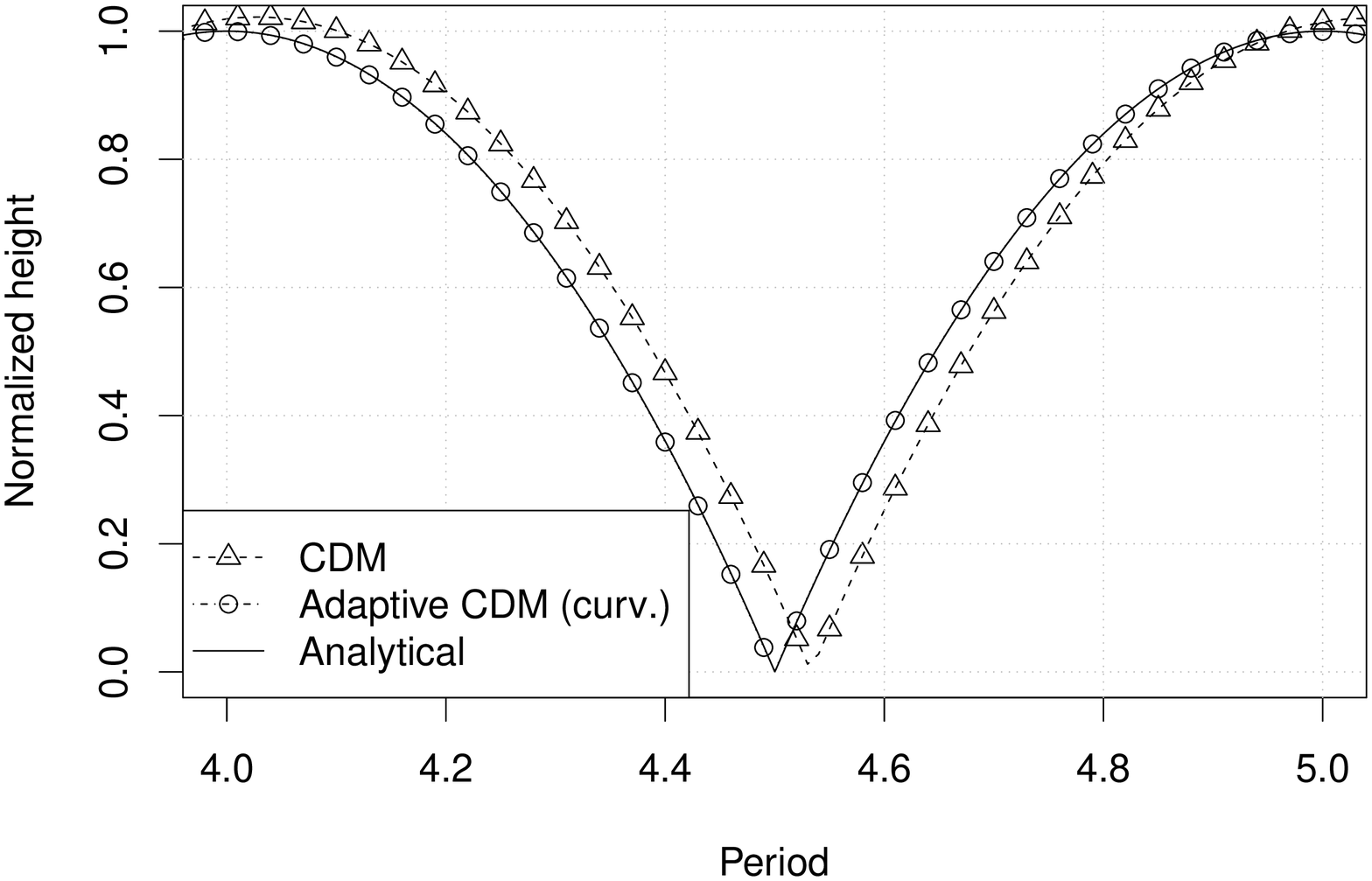}}
\subfigure[\label{fig:late}]{\includegraphics[width=0.45\linewidth]{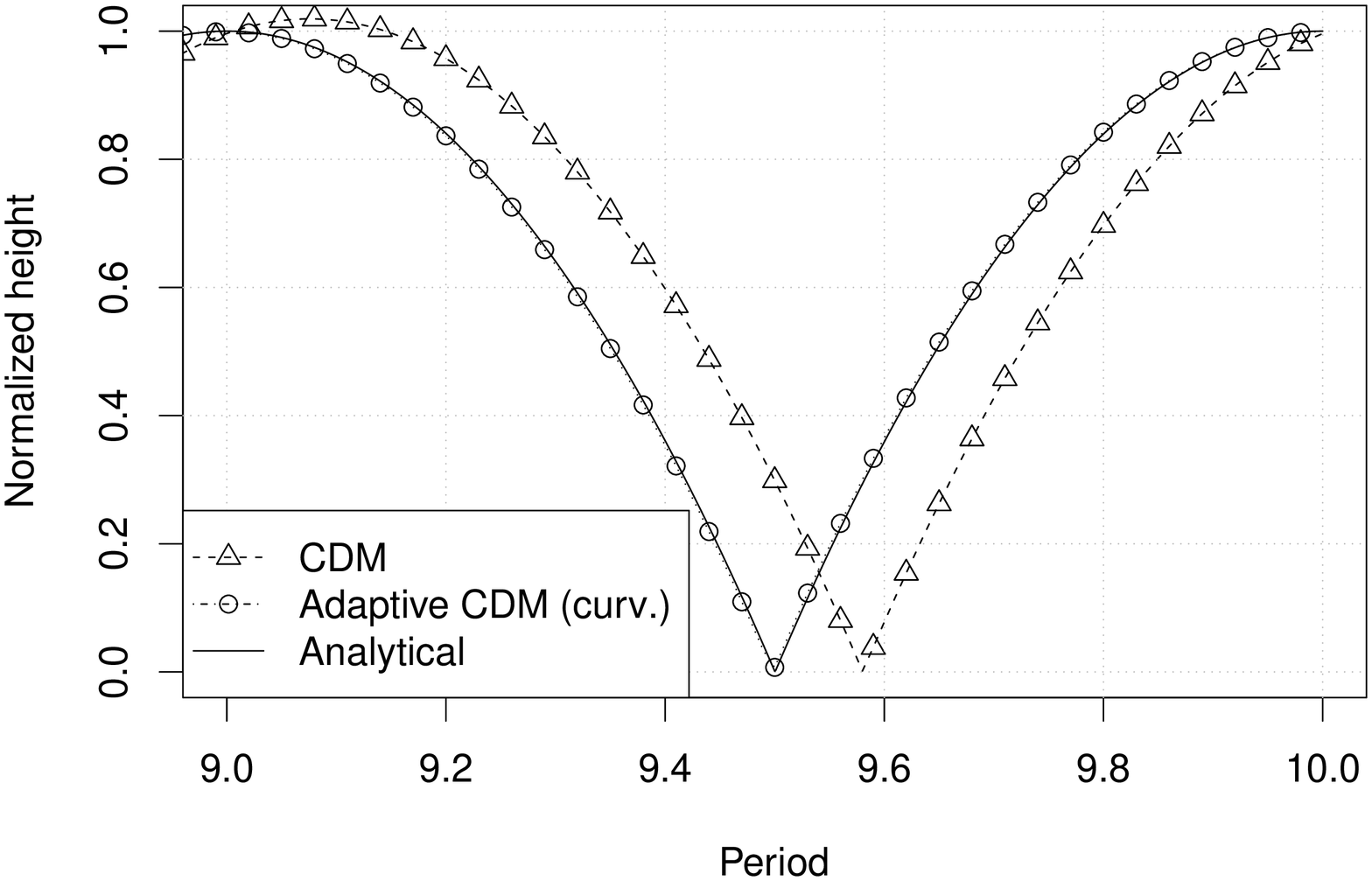}}
\caption{Displacement history.}
\label{fig:DEMmovement}
\end{figure}

Figure~\ref{fig:DEMmovement} presents the kinematics obtained from the numerical analyses, along with the analytic (exact) solution. 
The displacement history of the sphere is normalized by the launching height value $h_0$, and the time axis is normalized by the repetition period of the movement $t_f$. 

In the initial instants of the analysis, both adaptive and standard numerical solutions coincide with the exact response, cf. Figure~\ref{fig:initial}.
This is because of the compatibility of the time integrator used and the displacement function in the early stages of the simulation.
However, after the first contact, numerical responses become approximations of the analytical solution.
Note in~\eqref{eq:htparticle} that, when the particle is in contact with the surface, the displacement function uses trigonometrical functions.
The approximation by Taylor's expansion used by the temporal integrator introduces errors in these instants.
The errors are proportional to the value of the time increment used in the temporal integration, and increases in every successive contact. 

The standard algorithm strategy (CDM) has no control of the integration errors, and its response continues diverging  along all the simulation, cf. Figures~\ref{fig:interm} and~\ref{fig:late}.
The adaptive solution, since it adjusts the time increment interactively, has better precision than the algorithm that employs prescribed increments, as presented in Figure~\ref{fig:DEMmovement}.

The information of curvature history for the described problem is compared to the local error estimator suggested by Zienkiewicz and Xie~\cite{ZIENKIEWICZ1991}.
This algorithm and our proposal are comparable, since they both rely on the Newmark type scheme.
The local error can be defined, taking as reference the Newmark algorithm, based on the difference between two distinct approximations with different precisions to the displacement function.
The adaptive strategy based on natural apparent frequencies~\cite{PARK1980, PARK1980p2} is not suitable for this problem since the apparent current frequency is not zero only while in contact with the surface.
Zienkiewicz and Xie~\cite{ZIENKIEWICZ1991} noticed a similar problem in the analysis of some dissipating processes of consolidation.
This issue can be easily alleviated by using a local error or a maximum step size specified by the user.

\begin{figure}[hbt]
\centering	
\includegraphics[width=0.6\linewidth]{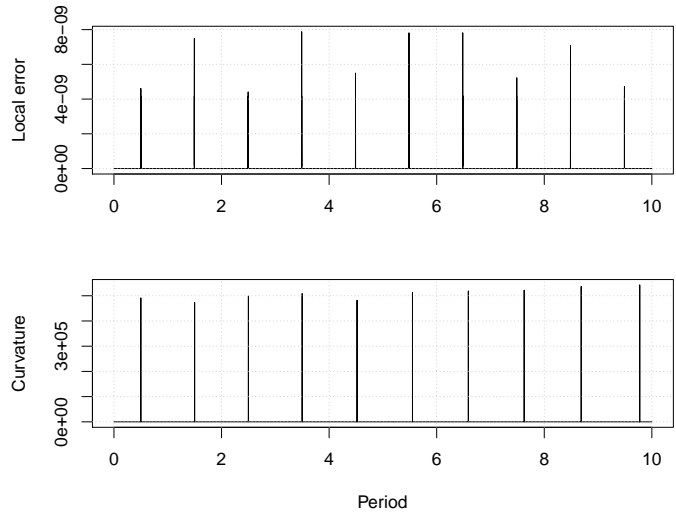}
\caption{Correlation: curvature \emph{versus} local error.}
\label{fig:localerrorandcurvaturecorrelation}
\end{figure}

Figure~\ref{fig:localerrorandcurvaturecorrelation} shows the displacement history curvature and the local error estimator along the time.
They are tightly related, since the peaks in these two functions occur at the same time.

\begin{figure}[hbt]
\centering
\includegraphics[width=0.6\linewidth]{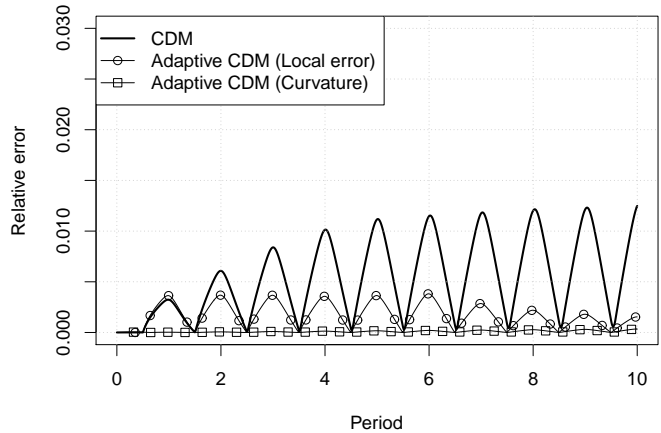}
\caption{Relative error history.}    
\label{fig:errorhistory}
\end{figure}

The relative errors of three integration strategies with respect to the exact solution is presented in Figure~\ref{fig:errorhistory}.
One of these strategies uses the standard algorithm without adaptation (CDM), and the others change the time step based on a refinement indicator: the local error~\cite{ZIENKIEWICZ1991} and the displacement history curvature. 
A close analysis of these curves reveals that all the strategies provide exact solutions until the first collision, and that errors appear after this event in all three solutions.
The peak error increases in the CDM and, for this particular simulation, decreases in the Adaptive CDM with local error estimator after a few periods.
Though barely visible, the error induced by the Adaptive CDM based on curvature also oscillates, but to an extent which is consistently at least two orders of magnitude smaller than the other two.

Figures~\ref{fig:localerrorandcurvaturecorrelation} and~\ref{fig:errorhistory} suggest that both curvature and local error reveal the instants when numerical integration is performed imprecisely, i.e., the moments when the sphere and the surface make contact.

\begin{figure}[hbt]
\centering
\includegraphics[width=0.6\linewidth]{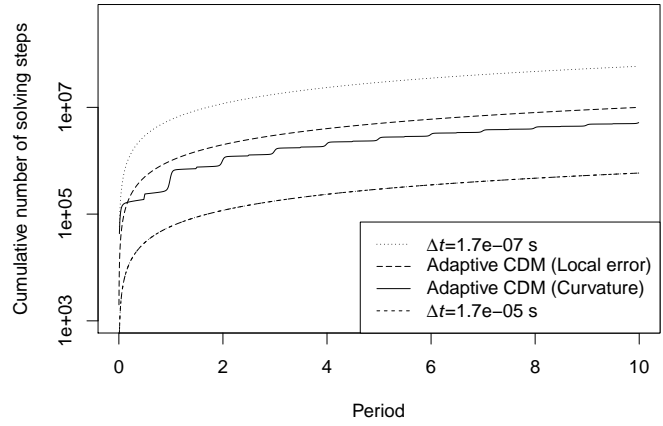}
\caption{Elastic collision: cumulative number of solving steps.}
\label{fig:dem_solvsteps_hist}
\end{figure}

In analogy with Figure~\ref{fig:fwd_solvsteps_hist}, Figure~\ref{fig:dem_solvsteps_hist} presents the cumulative number of solving steps of the four strategies considered in this example.
Again, the most and least demanding algorithms are, respectively, the ones with fixed time step and smallest and biggest increment, while the two adaptive solutions lie between those extremes.
After the initial stages of the simulation, the method based on curvature requires consistently less steps to achieve the solution than the strategy based on local error.
It is noteworthy that our proposal also produces better results, as presented in Figure~\ref{fig:errorhistory}.

\subsection{Varying the Integration Strategy}\label{subseq:VaryingIntegrationStrategy}

In this Section we compare the behavior of the curvature-based strategy for the same two problems above with three integration algorithms: CDM, EG-$\alpha$~\cite{Hulbert1996}, and Chung-Lee~\cite{Chung1994}.

\begin{figure}[hbt]
\centering
\subfigure[Four-wheel-dolly\label{fig:fwd_algorithms_comparing}]{\includegraphics[width=0.5\linewidth]{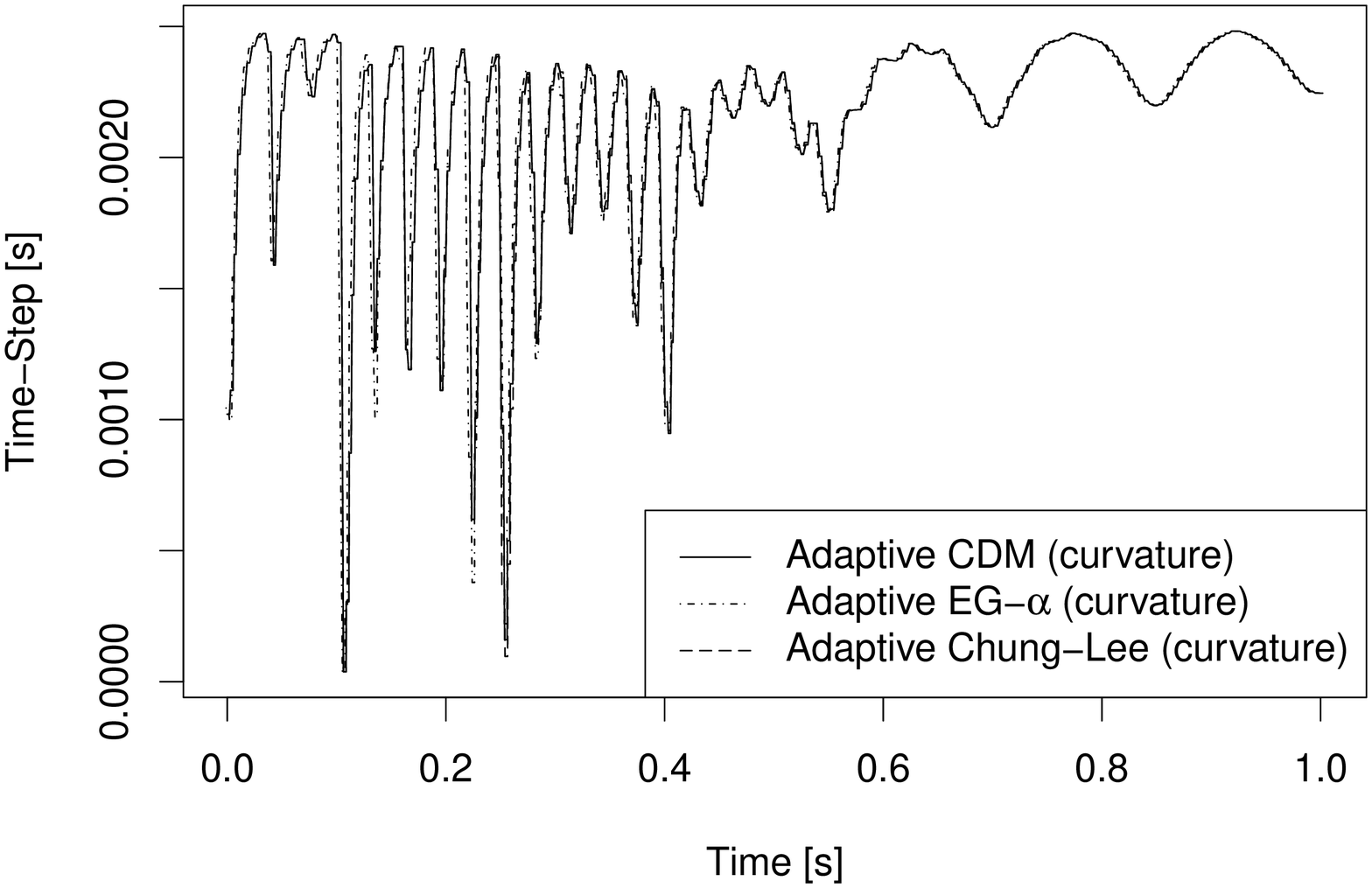}}
\subfigure[Elastic collision\label{fig:dem_algorithms_comparing}]{\includegraphics[width=0.5\linewidth]{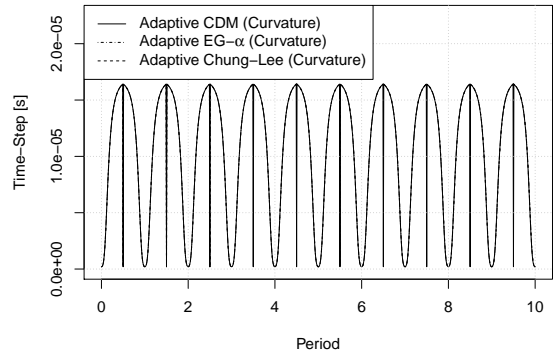}}
\caption{Time-steps }\label{fig:ComparingAlgorithms}
\end{figure}

Figures~\ref{fig:fwd_algorithms_comparing} and~\ref{fig:dem_algorithms_comparing} present the time steps required for solving the four-wheel dolly and elastic collision problems, respectively, using the aforementioned integration strategies.

Figure~\ref{fig:ComparingAlgorithms} shows that the time steps do not present significant changes when the integration strategy changes.
The strategy based on curvature captures the kinematic nature of the problem, rather than the particular technique used to solve it.
This experiment shows the versatility of our proposal.

\section{APPLICATION TO MOORING LINES}\label{sec:mooring}

The structural analysis of mooring lines has been target of intense research, mainly motivated by the oil industry.
This activity requires the use of floating offshore structures subjected to several environmental conditions as winds, waves and ocean currents.
The mooring lines are largely responsible for the stability of the vessel; they absorb loadings acting upon the submersible structure and transmit them to the marine soil, in some cases located thousands of meters from the ocean surface.

Studying the dynamics of such structures is the objective of software like \verb|Dynasim|~\cite{DYNASIM2002}, which employs a software component called \verb|DOOLINES|~\cite{Doolines2011}, an object oriented framework, developed in the C++ language, responsible for the dynamic analysis of mooring lines and risers.
This framework uses the Finite Element Method to make a spatial discretization of the lines, and Direct Integration Methods to solve the dynamics of the problem.
\verb|DOOLINES| employs the EG-$\alpha$ method.

\begin{figure}[hbt]
\centering
\includegraphics[width=0.6\linewidth]{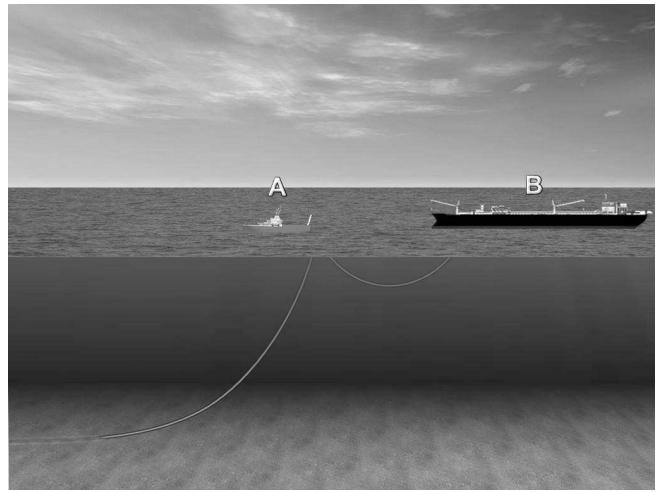}
\caption{Hookup scheme.}
\label{fig:hookupscheme}
\end{figure}

In the following example, the time-adaptive technique based on curvature is employed for mooring lines analysis of an installation procedure termed ``hookup'', using the values $b=10$, $\zeta=100$ and $\Delta t_{\max}=\unit[0.030342]{s}$.
In this procedure, a line portion initially suspended between ships  ``A'' and ``B'' has one of its terminations suddenly released by ship ``A'', as in Figure~\ref{fig:hookupscheme}. 

From the numerical modeling point of view, the release of the line corresponds to an instant of sudden change in the boundary conditions of the problem.
In general, sudden changes like this make an important contribution to the error in any numerical discretization.
The usual solution adopted is decreasing the increment of the integration ($\Delta t$).

In non adaptive algorithms, this decreasing is global and valid for the whole simulation.
This practice is simpler, but has negative impact in the time required to process the simulation, since it requires a greater number of solution steps.
Moreover, in many moments of the simulation, a larger time increment can be used without significant loss of the response quality.
Our adaptive algorithm acts identifying those critical instants, and suggesting according changes to the time increment.
It, therefore, improves the relation between quality and computational performance.

The adaptive algorithm based on curvature is employed in the following example that simulates the hookup procedure.
The objective is to investigate aspects related to the computational performance and numerical response quality.
The obtained data for the adaptive algorithm is compared to analyses where the time increment is fixed along the entire simulation.
For all analyses, a unique uniform spatial mesh is used to discretize the domain; it is composed of $\unit[5]{m}$ truss elements.
The physical and geometric parameters employed in the simulation are presented in Table~\ref{tab:hookupparameters}.

\begin{table}[hbt]
\caption{Physical and geometric parameters of the hookup model. }
\begin{center}
\begin{small}
\begin{tabular}{ll}
    \toprule
    Parameter &  Value\\
    \midrule
	Distance between ships			 & \unit[500]{m}\\
	Sea water level		                 & \unit[1000]{m}\\
	Suspended line length			 & \unit[550]{m}\\
	Total length of the line                 & \unit[1620]{m}\\
	Longitudinal stiffness of the line (EA)	 & \unit[1000]{kN}\\
	Weight in air (line)			 & \unitfrac[1.50]{kN}{m}\\
	Submerged weight (line)			 & \unitfrac[1.26]{kN}{m}\\
	Morison's drag coefficient 		 & 1.2\\
	Morison's inertia coefficient		 & 0.06\\
	Static friction coefficient 		 & 0.50\\
        Dynamical friction coefficient		 & 0.05\\
    \bottomrule
\end{tabular}
\end{small}
\end{center}
\label{tab:hookupparameters}
\end{table}

\begin{figure}[hbt]
\centering
\includegraphics[width=0.6\linewidth]{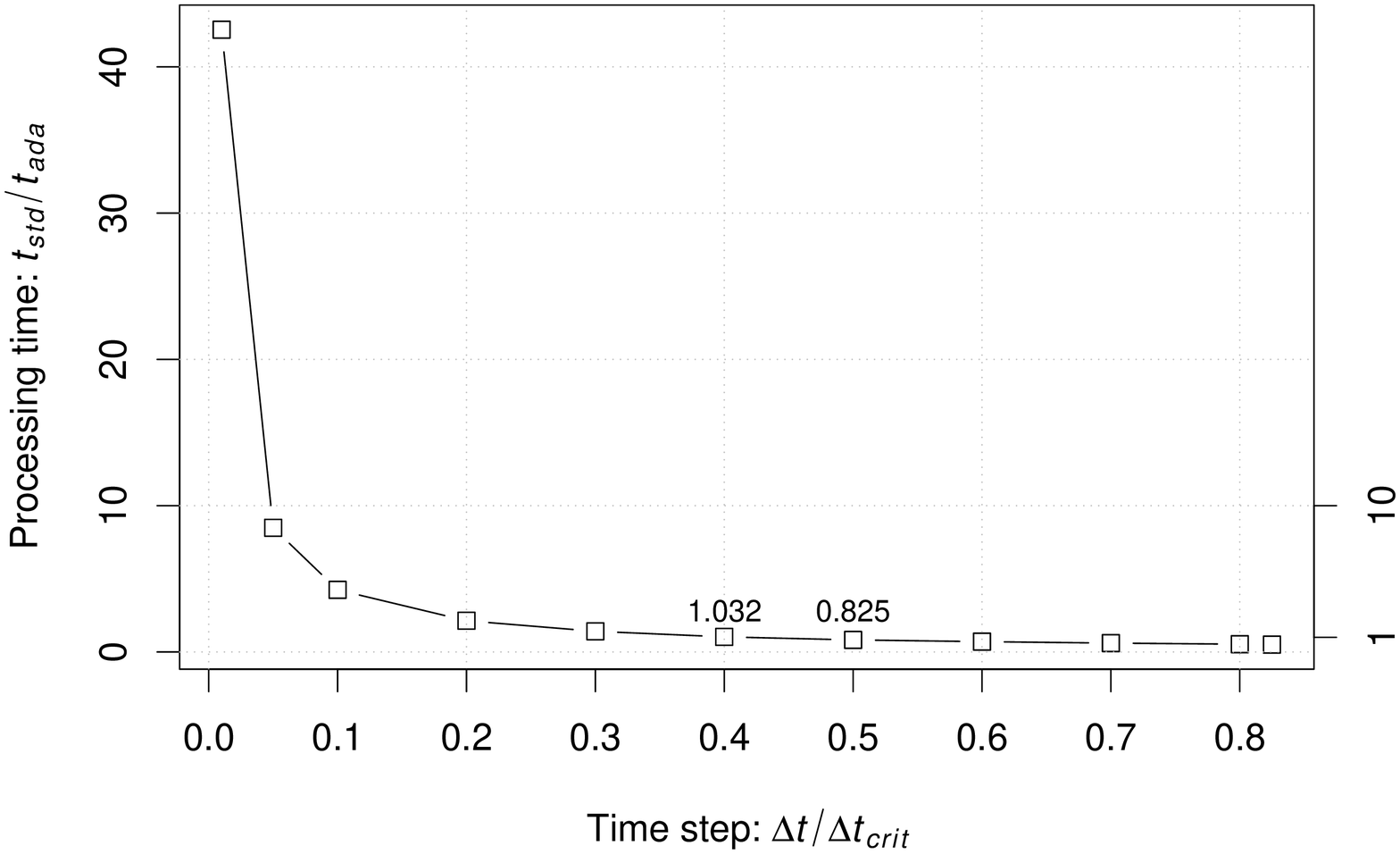}
\caption{Processing time comparison: standard \emph{versus} adaptive.}
\label{fig:hookupproctime}
\end{figure}

Figure~\ref{fig:hookupproctime} shows the computational time spent by the standard algorithm of time integration solving the hookup simulation ($t_{std}$) scaled by the time spent by the adaptive algorithm based in curvature ($t_{ada}$).
Additionally, the temporal discretization is shown proportional to the critical integration increment (${\Delta t}_{crit}$).
As expected, the processing time is inversely proportional to the time step used.
Besides, the adaptive algorithm has comparable performance to the standard algorithm with uniform mesh of size  $\Delta t = 0.4{\Delta t}_{crit}$.
In order to verify this assertion, two points were highlighted: $(0.4,1.032)$ and $(0.5,0.825)$.
These contiguous points are above and below the reference time; in most cases, the proposed algorithm performs better, in terms of computational time, than the standard one, and when it is outperformed by the latter, the difference in time is small.
We will see in the following that in those situations where the standard algorithm is faster than the adaptive one, the error produced by the former is bigger.
The experience was performed on an Intel Core 2 Duo@\unit[2.66]{GHz} processor, \unit[1]{GB} of RAM,  GCC~4.2.3 compiler and 32 bits architecture.

\begin{figure}[hbt]
\centering
\includegraphics[width=0.6\linewidth]{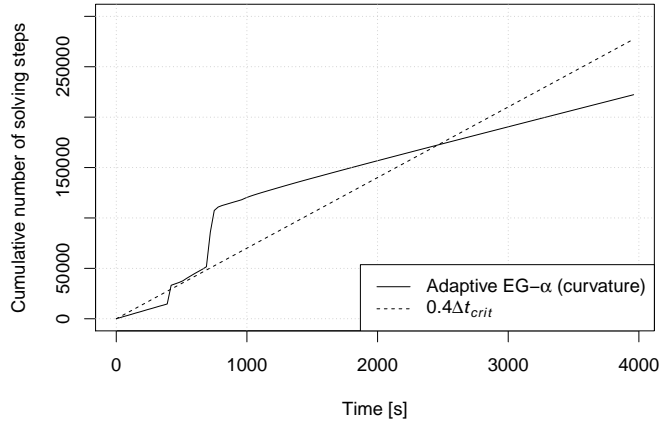}
\caption{Steps count comparison.}
\label{fig:stepscount}
\end{figure}

Figure~\ref{fig:stepscount} presents the number of steps that both the adaptive and the standard algorithms require for solving the elastic collision problem. 
By `step' we mean a call to the procedure that computes the external and internal forces.
The standard algorithm uses a fixed time step of $\Delta t = 0.4\Delta t_{crit}$. 

As presented in Figure~\ref{fig:stepscount}, the standard algorithm uses a number of solving steps which is linear with the time, whereas the adaptive algorithm begins demanding more operations and soon stabilizes in a more economic count.

\begin{figure}[hbt]
\centering
\includegraphics[width=0.6\linewidth]{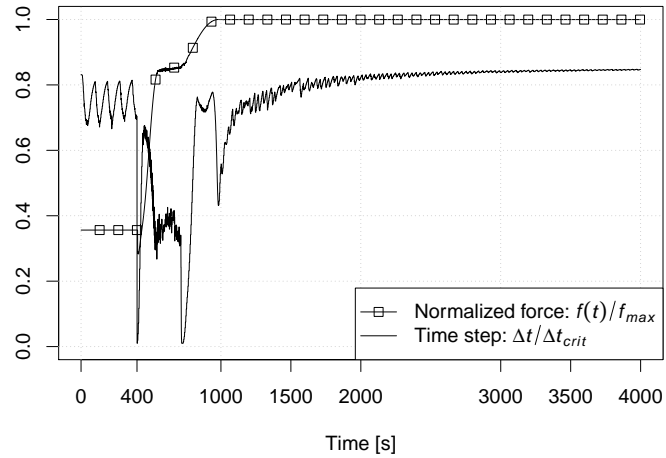}
\caption{Correlation between time increment and the behavior of the force in the top of line.}
\label{fig:hookupincrementhistory}
\end{figure}

Figure~\ref{fig:hookupincrementhistory} shows the time increment history for the adaptive algorithm.
In order to associate the changes of time increment with the dynamical behavior of the problem, the force acting at the end of the line (ship ``B'') is also presented.
The force is normalized by the maximum value $f_{\max}$, which corresponds to the one verified at the end of the analysis.

In the initial instants of the simulation, i.e., from $\unit[0]{s}$ to $\unit[400]{s}$, when the line is still hanging, the value of the force at the top of the line (Figure~\ref{fig:hookupscheme}, ship ``B'') is constant and defined only by the contribution of the weight of the suspended portion of the line.
When the line is released from ship ``A'' ($t=\unit[400]{s}$), the force has a small decrease because of the line motion, but then it increases until reaching its final value (approximately at \unit[$1000$]{s}), when the installation is done.

\begin{figure}[hbt]
\centering
\includegraphics[width=0.6\linewidth]{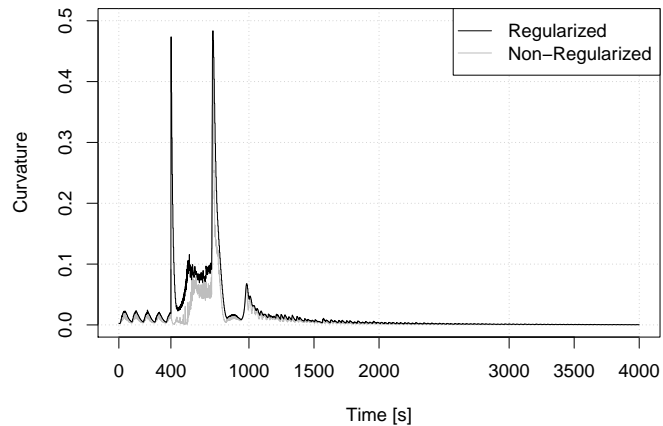}
\caption{Curvature history to the hookup problem.}
\label{fig:hookupcurvaturehistory}
\end{figure}

In terms of numerical errors, the instant when the line is released ($\unit[400]{s}$) and period of time when it collides with the marine soil ($\unit[400]{s}< t \leq \unit[1500]{s}$) are crucial in the analysis because of the changes of the line kinematics, as shown in  Figure~\ref{fig:hookupcurvaturehistory}.

Figure~\ref{fig:hookupcurvaturehistory} shows the displacement history curvature for both regularized and non-regularized methods.
The biggest values of curvature occur after the launching ($t>\unit[400]{s}$).
Note in Figure~\ref{fig:hookupincrementhistory} that the adaptive time increment is reduced at $\unit[400]{s}$, and kept smaller than the initial value, albeit oscillating, until stabilization at, approximately, $\unit[1500]{s}$.
These changes occur as a function of the variation in the displacement history.

\begin{figure}[hbt]
\centering
\includegraphics[width=0.6\linewidth]{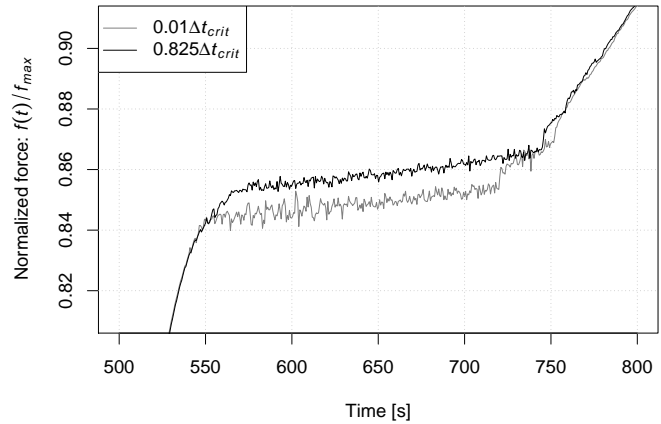}
\caption{Normalized forces obtained with $\Delta t = 0.825\Delta t_{crit}$ and with $\Delta t = 0.1\Delta t_{crit}$.}
\label{fig:hookuphisterrorstda}
\end{figure}

In the final stage of the simulation, when the line behavior is almost static, the adaptive algorithm increases the time increment value.
This value exceeds $0.80\Delta t_{crit}$ in some instants ($t>\unit[3500]{s}$ in Figure~\ref{fig:hookupincrementhistory}).
If a constant increment of $0.80\Delta t_{crit}$ was employed during all the simulation, the discretization error would be unacceptable in critical instants, as presented in Figure~\ref{fig:hookuphisterrorstda}.

Figure~\ref{fig:hookuphisterrorstda} presents the normalized forces, as produced by a poor time discretization ($\Delta t = 0.825\Delta t_{crit}$) and a refined one ($\Delta t = 0.01\Delta t_{crit}$) for the instants when the line collides with the soil.
The refined solution is close to the analytical one, so the difference between the curves verifies that integration errors contaminate the solution obtained using $\Delta t = 0.825\Delta t_{crit}$.

\begin{figure}[hbt]
\centering
\includegraphics[width=0.6\linewidth]{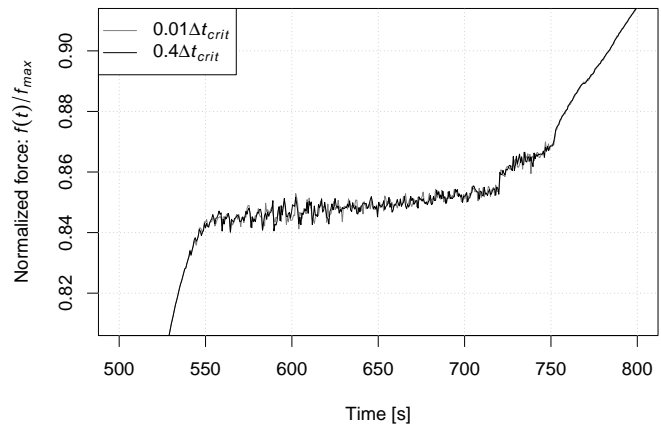}
\caption{Normalized forces obtained with $\Delta t = 0.01\Delta t_{crit}$ and with $\Delta t = 0.4\Delta t_{crit}$.}
\label{fig:hookuphisterrostdb}
\end{figure}

For some design situations, even the presented poor discretization may yield good enough results.
However, more refined and, therefore, dependable results are often necessary.
A solution for this issue would be employing reduced increments, as commented before and illustrated in Figure~\ref{fig:hookuphisterrostdb}, where $\Delta t = 0.40\Delta t_{crit}$ was employed.

Note in Figure~\ref{fig:hookuphisterrostdb} that the responses produced by the standard algorithm using $\Delta t = 0.01\Delta t_{crit}$ and $\Delta t = 0.40\Delta t_{crit}$ are alike, suggesting that the latter is enough to obtain a response of good quality.
Choosing the time increment that ensures a good precision of response however is not an easy task, as previously shown.

\begin{figure}[hbt]
\centering
\includegraphics[width=0.6\linewidth]{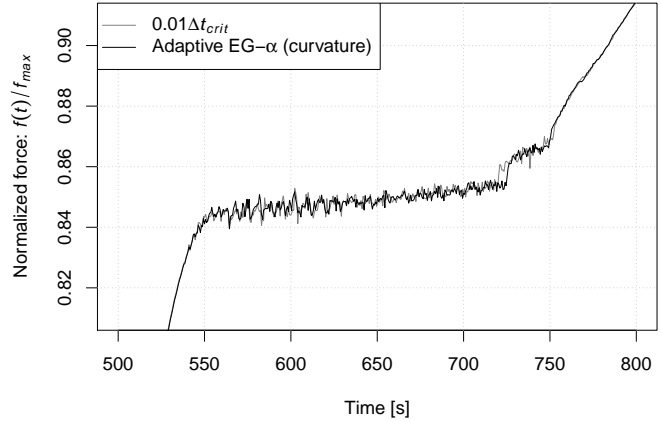}
\caption{Normalized forces obtained with $\Delta t = 0.01\Delta t_{crit}$ and with adaptive algorithm.}
\label{fig:hookuphisterrorada}
\end{figure}

Figure~\ref{fig:hookuphisterrorada} presents a comparison of the force response at the top of the line, as obtained by curvature based and standard algorithms.
The latter used the same refined discretization of reference ($\Delta t = 0.01\Delta t_{crit}$).
Note that the difference of the responses is small, suggesting that the information provided by the geometric indicator can be successfully employed to improve the quality of the response by controlling the time step.
Besides the automatic improvement in quality verified in Figure~\ref{fig:hookuphisterrorada}, the processing time required by the adaptive algorithm noticeably smaller (Figure~\ref{fig:hookupproctime}).
This association between quality and computational performance is an important requirement for any efficient adaptive technique of numerical time integration.

\section{CONCLUSIONS}\label{sec:conclu}

The time-adaptive procedure presented here can be used in the solution of dynamics problems that employ iterative methods, as the Direct Integration Method.
As presented, the refinement estimator requires only a few additional computations, and it is able to identify kinematic changes which are usually related to errors in temporal numerical integration.
Our proposal allows an adaptive control of the time step of integration and, therefore, improves the relation between quality of response and time processing.
We observed a good association between the geometric indicator curvature and the error estimator based on approximate local solutions.

The adaptive procedure proposed is appropriate for dynamics problems which require localized refinement for the temporal meshes.
Simulations involving abrupt changes in boundary conditions or in the unknown kinematics are good examples of possible successful applications of the adaptive technique.

This article presents results which were obtained following the Reproducible Research guidelines~\cite{RRSignalProcessing}.
All the relevant information is available at \url{http://loi.lccv.ufal.br/art-ada}.

\acks

The authors would like to thank CNPq, CAPES and FAPEAL that provided partial support for this work.

\bibliographystyle{wileyj}
\bibliography{references}

\end{document}